\documentclass[11pt,twoside,english,reqno,a4paper]{amsart}

\usepackage{amssymb, amsmath}
\usepackage{color}

\definecolor{red}{rgb}{0.9,0,0}
\usepackage{palatino}
\definecolor{azul}{rgb}{0,0,1}

\usepackage[colorlinks=false]{hyperref}
\usepackage{tikz-cd}
\usetikzlibrary{patterns}

 \makeatletter
\@namedef{subjclassname@2020}{%
  \textup{2020} Mathematics Subject Classification}
\makeatother

\usepackage{indentfirst}
\usepackage{graphicx}

\newtheorem{theorem}{Theorem}[section]
\newtheorem{proposition}{Proposition}[section]

\newtheorem{lemma}{Lemma}[section]
\newtheorem{remark}{Remark}[section]
\newtheorem{definition}{Definition}[section]
\newtheorem{example}{Example}[section]
\usepackage[left=2.9cm,right=2.9cm,top=2.8cm,bottom=2.8cm]{geometry}

\usepackage{geometry}
\def\re{\mathbb{R}}
\def\bk{\color{black}}
\def\dem{\noindent{\it Proof. }}

\newcommand{\fim}{\hfill $\Box$}

\usepackage{booktabs,amsmath}
\def\Xint#1{\mathchoice
    {\XXint\displaystyle\textstyle{#1}}%
    {\XXint\textstyle\scriptstyle{#1}}%
    {\XXint\scriptstyle\scriptscriptstyle{#1}}%
    {\XXint\scriptscriptstyle\scriptscriptstyle{#1}}%
      \!\int}
\def\XXint#1#2#3{{\setbox0=\hbox{$#1{#2#3}{\int}$}
    \vcenter{\hbox{$#2#3$}}\kern-.5\wd0}}
\def\dashint{\Xint-}

\title[Hardy potentials as $p\to1^+$: the critical drift case]{Existence, non-existence and degeneracy of limit solutions to  $p$-Laplace problems involving Hardy potentials as $p\to1^+$. \\  The case of a critical drift}

\author{J. C. Ortiz Chata }
\address{Departamento de Matem\'atica, Universidade Federal de São Carlos - UFSCar, 13565-905, S\~ao Carlos - SP, Brazil}
\address{juanchata@ufscar.br}
\author[F. Petitta]{Francesco Petitta}
\address{Dipartimento di Scienze di Base e Applicate per l' Ingegneria, Sapienza Universit\`a di Roma, Via Scarpa 16, 00161 Roma, Italia}
\address{francesco.petitta@uniroma1.it}

\date{}
\begin{document}
\pretolerance10000

\maketitle

\begin{abstract}
In this paper we analyze the asymptotic behaviour as $p\to 1^+$ of solutions $u_p$ to
$$
\left\{
\begin{array}{rclr}
-\Delta_pu&=&\lambda|\nabla u|^{p-2}\nabla u\cdot\frac{x}{|x|^2}+ f&\quad \mbox{ in } \Omega,\\
u_p&=&0 &\quad \mbox{ on }\partial\Omega,
\end{array}\right.
$$
where $\Omega$ is a bounded open subset of $\re^N$ with Lipschitz boundary containing the origin, $\lambda\in\re$, and $f$ is a nonnegative datum in $L^{N,\infty}(\Omega)$. 
As a consequence, under suitable smallness assumptions on $f$ and $\lambda$, we show sharp  existence results  of bounded solutions to the Dirichlet problems 
$$\begin{cases}
					 \displaystyle - \Delta_{1} u = \lambda\frac{D u}{|D u|}\cdot \frac{x}{|x|^2}+f &  \text{in}\, \Omega, \\
					 u=0 & \text{on}\ \partial \Omega,			
				 \end{cases}
$$
 where $\displaystyle \Delta_{1}u=\hbox{div}\,\left(\frac{Du}{|Du|}\right) $ is the $1$-Laplacian operator. The case of a generic drift term in $L^{N,\infty}(\Omega)$ is also considered. 
   Explicits examples are given in order to show the optimality of the main assumptions on the data.
 \end{abstract}

\numberwithin{equation}{section}
\bibliographystyle{plain}
\maketitle
\tableofcontents

\section{Introduction}
Following \cite{OrtizPetitta2024},  the present  paper  represents the second of a series of studies concerning $1$-Laplace type problems in presence of Hardy type nonlinear potential terms.  The main motivation of this paper concerns the understanding of problems that reads as  
\begin{equation}
\label{Pintro} 
\left\{
\begin{array}{rclc}
\displaystyle-\Delta_1u&=&\lambda\frac{D u}{|D u|}\cdot \frac{x}{|x|^2}+f&\quad \mbox{ in } \Omega,\\
u&=&0 &\quad \mbox{ on }\partial\Omega,
\end{array}
\right.
\end{equation} 
where $\Delta_1u=\hbox{div}\,\left(\frac{Du}{|Du|}\right)$, $\Omega\subset\mathbb{R}^N$ $(N\geq 2)$ is a open set containing the origin with bounded Lipschitz boundary, $0<|\lambda|<N-1$, and $f$ is in $L^{N,\infty}(\Omega)$  with suitable small norm. In the slightly less general, but more transparent,  case of a datum  $f\in L^{N}(\Omega)$ this smallness condition reads as 
\begin{equation}
\label{Hintro}
{S_N}\|f\|_{L^N(\Omega)}+\frac{|\lambda|}{N-1}\leq1,
\end{equation}    
where $S_N$ is the best constant in Sobolev inequality (see \eqref{sn} below).

\smallskip 
Notice that the drift term $\frac{Du}{|D u|}\cdot \frac{x}{|x|^2}$ is critical in many senses. First of all the ratio  $\frac{Du}{|Du|}$ (that appears twice in the equation) may have  no sense at those points where $D u =0$; the situation is even worst as,  in general, the natural functional space in which these type of problems are built-in is the space of function with bounded variation and so $Du$ may be only a Radon measure vector. The question on how to interpret the singular quotient  $\frac{Du}{|Du|}$ is nowadays well understood since the works of \cite{Demengel1999, acm2001, AndreuMazonMollCaselles2004, AndreuCasellesMazon2004} in which $\frac{Du}{|Du|}$ is formally replaced by a vector field  ${\bf z}\in L^N(\Omega, \mathbb{R}^N)$ that plays its role, using Anzellotti's paring theory \cite{Anzellotti1983}.

 On the other hand,   $\frac{x}{|x|^2}$ does   not belong to $L^N(\Omega,\mathbb{R}^N )$, but only to the Marcinkiewicz (or Lorentz) space $L^{N,\infty}(\Omega,\mathbb{R}^N)$ and, since as we said  $BV(\Omega)$ is the natural space in which problem \eqref{Pintro} is set, it becomes a borderline case   that  needs to be handled by mean of Hardy's inequality (see \eqref{HIe} below). 
 
 As we mentioned problem \eqref{Pintro} is in some sense the first order counterpart of   problems with critical Hardy zero order terms  
\begin{equation}
\label{op24} 
\left\{
\begin{array}{rclr}
\displaystyle-\Delta_1u&=&\frac{\lambda}{|x|}\frac{u}{|u|}+f&\quad \mbox{ in } \Omega,\\
u&=&0 &\quad \mbox{ on }\partial\Omega,
\end{array}
\right.
\end{equation} 
studied in \cite{OrtizPetitta2024}, since both perturbation terms share the same singular Hardy type  critical potential growth $\frac{1}{|x|}$. 

\smallskip 
Problems involving the  $1$-Laplace operators are widley  employed in a huge number of applications as  image restoration since Rudin, Osher and Fatemi in \cite{ROF} whose goal consisted in proposing a model for  the reconstruction  of an image $u$ given a blurry one $u_0=Au+n$; here A is a linear operator, e.g. a blur,  and $n$ is a random noise.   This problem consists in  looking for minimizers of functionals as  
\begin{equation}\label{rof}
\int_{\Omega} |\nabla u| dx + \int_{\Omega}|A u- u_0|^2 dx,
\end{equation}
which are naturally relaxed in $BV(\Omega)$ in order to reinforce
 edges. Formally the sub-differential of the relaxed energy term in \eqref{rof} coincide with the $1$-Laplace operator.  

Other possible applications of $1$-Laplace type problems  include the study of  torsional creep of  a cylindrical bar of constant cross section in $\re^2$, as well as  more theoretical issues of geometrical nature;  an  account on these and more possible applications can be found in  \cite{OsSe, Kawohl1990, Kawohl1991},   \cite{Sapiro, M, BCRS},  the monograph \cite{AndreuCasellesMazon2004}, and references therein. 

From the mathematical point of view, problems as \eqref{Pintro} do not have in general a valid variational structure so that direct methods of calculus of variations  are not directly applicable.  One of the ways to overcome this difficulty, it is to study the solution through  ``approximate" $p$-Laplacian problem solutions as $p$ go to $1^+$, as in \cite{Kawohl1990, ct, mst, OrtizPetitta2024}.

The case $\lambda=0$ in \eqref{Pintro} has been extensively addressed. 
The class of problems as  in \eqref{Pintro} when $\lambda=0$  have been studied in a series of papers \cite{Kawohl1990, ct, mst} as an outcome of the study the behaviour of the solutions to the problem
\begin{equation}
\label{PLa}
\left\{
\begin{array}{rclr}
-\Delta_pu&=&f&\quad \mbox{ in } \Omega,\\
u&=&0 &\quad \mbox{ on }\partial\Omega,
\end{array}
\right.
\end{equation} 
as $p\to 1^+$ whenever the norm of $f$ is small. For instance, in \cite{Kawohl1990},  the author studied the existence of variational solutions for case $f=1$ provided $\Omega$ is suitably small. In particular, if \eqref{PLa} is formulated as a variational problem it is shown that it could admit  a non-trivial minimizer. In \cite{ct} the authors prove    the degeneracy $u_p\to 0$ as $p\to1^{+}$ provided  the norm of $f$ is too small, while in  \cite{mst} it is shown that the  limit of $u_p$ no matters of its degeneracy satisfies the concept of solution introduced in  \cite{acm2001, AndreuMazonMollCaselles2004, AndreuCasellesMazon2004}; in the same paper also  (non-trivial) solutions are obtained for data satisfying  $\|f\|_{W^{-1,\infty}(\Omega)}=1$. The case of problems as in \eqref{PLa} with   $L^1$ data may present very degenerate situations  (\cite{MST2009}) unless   in presence of regularizing  lower order terms (see for instance \cite{lops,bop}).  

Concerning problem \eqref{op24} an  optimal description of solutions depending on the size of both the parameter $\lambda$ and the datum $f$,  e.g.  \eqref{Hintro} in case $f\in L^{N}(\Omega)$, by mean of the study of the asymptotic behaviour of $p$-Laplace associated problems, and inspired by \cite{AzoreroAlonso1998, BoccardoOrsinaPeral2006} for the case $p>1$,  is the content of \cite{OrtizPetitta2024}. 

The case of a negative parameter $\lambda$ is of particular interest also for its connection with singularly weighted elliptic problems  as 
\begin{equation}
\label{Pa}
\left\{
\begin{array}{rclc}
-{|x|^a}\hbox{div}\,\left(\frac{1}{|x|^a}\frac{Du}{|Du|}\right)&=&f&\hbox{ in }\; \Omega\\
              u&=&0&\hbox{ on  }\;\partial\Omega,
\end{array}
\right.
\end{equation}
where  $0<a<N-1$ which are naturally set (here $a=-\lambda$), via Caffarelli-Kohn-Nirenberg type's inequality in weighted Sobolev spaces (see \cite{OrtizPimentaSegura2021} and Section \ref{sbs6} below). 

Finally, problems involving $1$-Laplace problems and gradient type lower order terms with natural growth have also been considered; without the aim to be complete we mention \cite{4,MS2013,gop,bop} (and references therein).

The goal of this work is to study, under minimal assumptions, the asymptotic behaviour of solutions to
\begin{equation}\label{rclc}
\left\{
\begin{array}{rclc}
-\Delta_pu&=&\lambda|\nabla u|^{p-2}\nabla u\cdot\frac{x}{|x|^2}+ f& \hbox{ in }\;\Omega,\\
              u&=&0&\hbox{ on }\;\partial\Omega,
\end{array}
\right.
\end{equation}
as $p\to 1^+$. As a result sharp existence, non-existence and regularity  results for solutions to \eqref{Pintro} are derived.  

 Existence and regularity of solutions of \eqref{rclc}  has been studied in  \cite{LeonoriPetitta2007} where under  smallness assumptions  on the parameter $\lambda$ (see also \cite{BMMP2003, DV1995, LMPP2011, BG1989} for related problems).

This paper is organised as follows. In Section \ref{sec2} we set the main notation and we give some preliminaries on  Lorentz spaces, bounded variations functions space, i.e.  $BV(\Omega)$,  and the extended Anzellotti-Chen-Frid  theory of pairings. In Section 3 we state  the main assumptions and results of the paper on, respectively existence of solutions for problem \eqref{Pintro}  in case of a datum $f$ in $L^N(\Omega)$ (Theorem \ref{exi}), asymptotic behaviour of the solutions of \eqref{rclc} (Theorem \ref{AsympB}), and boundedness of the obtained solutions (Theorem \ref{bo}) provided $f$ in $L^q(\Omega)$, $q>N$. Section \ref{mains} is devoted to the proofs of those results. Later on we show how the presented results can be extended, respectively,  to the case of a source $f$ in $L^{N, \infty}(\Omega)$ and $W^{-1,\infty}(\Omega)$  (Section \ref{sec5}) and to the case of a generic drift term (Section \ref{sec6}). Finally, in   Section \ref{sec7}, also through explicit examples,  we investigate the optimality of our assumptions along with some  non-existence results for large data.

\section{Preliminaries}\label{sec2}
\noindent{\bf Notations.}

We   denote by
\begin{itemize}
\item $\mathcal{H}^{N-1}$, the $N-1-$dimensional Hausdorff measure on $\mathbb{R}^N$;
\item $ \mathcal{L}^N$, the $N-$dimensional Lebesgue measure on $\mathbb{R}^N$ and for each measurable set $E\subset \mathbb{R}^N$, we write $|E|$ instead of $\mathcal{L}^N(E)$;
\item $B_r(x)$, the ball in $\mathbb{R}^N$ with center in $x$ and radius $r$;
\item $C_N$, the $N-$dimensional Lebesgue measure of $B_1(0)$, i.e., $C_N=|B_1(0)|$;
\item $C_c(\Omega)$, the set of continuous functions with compact support in $\Omega$;
\item $C_0(\Omega)$, the completion of $C_c(\Omega)$ with respect to $\|u\|=\sup\{|u(x)|\colon x\in \Omega\}$.

\end{itemize}
Also, for $k\geq 1$, we denote  by $T_k\colon \mathbb{R}\to \mathbb{R}$ the truncation function  defined by 
\begin{equation*}
T_k(s)=\left\{
\begin{array}{lr}
s&\hbox{ if }\; |s|\leq k\\
k\frac{s}{|s|}&\hbox{ if }\; |s|>k.
\end{array}
\right.
\end{equation*}
We  denote  by $\hbox{Sgn}\colon \mathbb{R}\to \mathcal{P}(\mathbb{R})$ the set-valued sign function whose definition is given by
\begin{equation*}
\hbox{Sgn}(s)=\left\{
\begin{array}{lr}
\frac{s}{|s|}& \hbox{ if }\; s\neq 0\\
 {[0,1] } & \hbox{ if }\; s=0.
\end{array}
\right.
\end{equation*}
\bk

Regarding the integrals notation, if there is no ambiguity we will use the notation
$$
\int_\Omega f:=\int_\Omega f(x)\,dx
$$
and, if $\mu$ is a Radon measure,
\[\int_\Omega f\mu:=\int_\Omega f\, d\mu\,.\]

\medskip

We underline the use of the standard convention to do not relabel an extracted compact subsequence.

\subsection{Lorentz spaces}

Let $u$ be a measurable function in an open bounded set $\Omega\subset \mathbb{R}^N$. We recall that the distribution function, the non-increasing rearrangement and the spherically non-increasing rearrangement of $u$ are given by 
\begin{equation*}
\label{ReaU}
\begin{array}{rclc}
\alpha_u(s)&=&|\{x\in \Omega\>: \> |u(x)|>s\}|&\forall\, s\geq 0;\\
u^*(t)&=&\sup\{s>0\>:\> \alpha_u(s)> t\}&\forall\;t\in (0, |\Omega|);\\
u^{\#}(x)&=&u^*(C_N|x|^N)&\forall\; x\in \Omega^{\#};
\end{array}
\end{equation*}
respectively, where $\Omega^{\#}$ is the ball centred at the origin having the same measure as $\Omega$.

The following result is proven in \cite[Theorem 378]{HLP1964}. 
{
\begin{proposition}
\label{Inq}
Let $u, v$ be nonnegative measurable functions on $\Omega$. Then the following inequality holds:
\begin{equation*}
\int_{\Omega}u(x)v(x)\,dx\leq \int_{0}^{\infty}u^*(t)v^*(t)\,dt.
\end{equation*}
\end{proposition}
}

\begin{definition}
Let $1\leq p <\infty$ and $1\leq q \leq \infty$.  The Lorentz space, denoted by $L^{p, q}(\Omega)$, is defined as 
$$
L^{p,q}(\Omega)=\{ u\>:\> u\mbox{ measurable on }\>\Omega,\; \|u\|_{L^{p,q}(\Omega)}<\infty\},
$$
where 
\begin{equation*}
\|u\|_{L^{p,q}(\Omega)}=\left\{
\begin{array}{lcr}
\left[\int_{0}^{\infty}[t^{1/p}u^*(t)]^q\frac{dt}{t}\right]^{1/q}& \hbox{ if } &1\leq p<\infty,\; 1\leq q<\infty\\
\sup_{t>0}t^{1/p}u^*(t)& \hbox{ if } & 1\leq p< \infty,\; q=\infty.
\end{array}
\right.
\end{equation*}
\end{definition}

\begin{proposition}[H\"older's inequality]
\label{HoI}
Let $1<p<\infty$, $ 1\leq q \leq  \infty$, $1/p+1/p'=1$ and $1/q+1/q'=1$. Let $u\in L^{p,q}(\Omega)$ and $v\in L^{p',q'}(\Omega)$. Then 
$$
\left|\int_{\Omega}u(x)v(x)\,dx\right|\leq \|u\|_{L^{p, q}(\Omega)}\|v\|_{L^{p', q'}(\Omega)}.
$$
\end{proposition}

The proof of the previous proposition can be found in \cite{ONeil1968} (see also \cite[Appendix]{DAnconaFanelli2007}).

\begin{proposition}
\label{SL}
Let $1\leq p <N$ and let $q\in [1, p]$. Then there exists a constant $\gamma_{p}>0$ such that the following inequality holds for all $u\in C^\infty_c(\mathbb{R}^N)$
$$
\|u\|_{L^{p^*,q}(\Omega)}\leq \gamma_p\|\nabla u\|_{L^{p,q}(\Omega)},
$$
where $ p^*=Np/(N-p)$. Futhermore, the constant
\begin{equation*}
\label{gp}
\gamma_p=\frac{p}{(N-p)C_N^{1/N}}
\end{equation*}
is optimal in this inequality.
\end{proposition}

This result is classical and its proof is in \cite{Alvino1977}.

\subsection{The space $BV$}

The space of functions of bounded variation over $\Omega$, denoted by $BV(\Omega)$, is made up of all functions in $L^1(\Omega)$ whose distributional derivative is a finite Radon measure in $\Omega$, i.e.,
$$
BV(\Omega)=\{u\in L^1(\Omega)\>;\> Du\in \mathcal{M}(\Omega, \mathbb{R}^N)\}.
$$

Equivalently, $u$ is a bounded variation function in $\Omega$ if and only if $u\in L^1(\Omega)$ and the following quantity
\begin{equation}
\label{TV}
 \sup\left\{\int_{\Omega}u\,\mbox{div}\,\varphi\,dx\>:\>\varphi\in C^1_c(\Omega,\mathbb{R}^N),\quad|\varphi|\leq 1\right\}
\end{equation}
is finite. Moreover, this quantity is equal to the total variation of $Du$ on $\Omega$, i.e., 
$$
\displaystyle |Du|(\Omega)=\int_{\Omega}|Du|=\sup\left\{\int_{\Omega}u\,\mbox{div}\,(\varphi)\,dx\>:\>\varphi\in C^1_c(\Omega,\mathbb{R}^N),\quad|\varphi|\leq 1\right\}.
$$
 For more details we refer to \cite[Proposition 3.6]{AmbrosioFuscoPallara}.

The space $BV(\Omega)$ is a Banach space endowed with the norm
$$
\|u\|_{BV(\Omega)}=\int_{\Omega}|Du|+\int_{\Omega}|u|\,dx.
$$

On the other hand, the classical result of existence of trace operator between Sobolev spaces $W^{1,p}(\Omega)$ and $L^1(\partial\Omega)$ can be extended between $BV(\Omega)$ and $L^1(\partial\Omega)$ (for instance see \cite{AmbrosioFuscoPallara, AttouchButtazzoMichaille2006, EvansGariepy2015}). The following holds:
\begin{proposition}
\label{BT}
Let $\Omega\subset \mathbb{R}^N$ be an open set with bounded Lipschitz boundary $\partial\Omega$ and $u\in BV(\Omega)$. Then, for $\mathcal{H}^{N-1}-$almost every $x\in \partial\Omega$ there exists $u^\Omega(x)\in \mathbb{R}$ such that
$$
\lim_{\rho\to 0}\rho^{-N}\int_{\Omega\cap B_\rho(x)}|u(y)-u^\Omega(x)|\,dy=0.
$$
Moreover, 
$$
\|u^\Omega\|_{L^1(\partial\Omega)}\leq C\|u\|_{BV(\Omega)},
$$
for some constant $C$ depending only on $\Omega$.
\end{proposition}

\begin{definition}[Strict convergence]
Let $u, u_n\in BV(\Omega)$, for all $n\in \mathbb{N}$, we say that $u_n$ strictly converges in $BV(\Omega)$ to $u$ if
\begin{eqnarray*}
u_n\to u\quad \hbox{ in }\; L^1(\Omega),\\
\int_{\Omega}|Du_n|\to \int_{\Omega}|Du|,
\end{eqnarray*}
as $n\to \infty$.
\end{definition}

It can be seen  that trace operator $u\mapsto u^\Omega$ with the topology induced by the strict convergence  is continuous on every open subset of $\mathbb{R}^N$ with bounded Lipschitz boundary. For simplicity from now on we will write $u$ ($\in L^1(\partial\Omega)$) instead of $u^\Omega$.

\begin{remark}
Since the trace operator $ u\mapsto u\in L^1(\partial\Omega)$ is continuous, using \cite[Proposition 2]{Miranda1974/75},  the above norm is equivalent to the following
$$
\|u\|:=\int_{\Omega}|Du|+\int_{\partial\Omega}|u|\,d\mathcal{H}^{N-1}.
$$
 \end{remark}

The space $C^\infty(\Omega)$ is dense in $BV(\Omega)$ with respect to strict convergence (see \cite[Theorem 3.9]{AmbrosioFuscoPallara}). Moreover one has the following  weak lower semicontinuity result that  immediately follows by the definition in \eqref{TV}.

\begin{proposition}
Let $u_n$ be a bounded sequence in $BV(\Omega)$ such that $u_n\to u$ in $L^1(\Omega)$. Then one has that $u\in BV(\Omega)$ and 
$$
\int_{\Omega}|Du|\leq\liminf_{n\to\infty}\int_{\Omega}|Du_n|.
$$
\end{proposition}

\begin{lemma}[Caffarelli-Kohn-Nirenberg inequality]
\label{CKN}
Let $p\geq 1$. Let $0\leq a<N-1$ and $0\leq b\leq 1$. Then there exists a constant $C>0$ such that 
$$
\left(\int_{\mathbb{R}^N}\frac{1}{|x|^{bp^*_b}}|\varphi|^{p^*_b}\,dx\right)^{1/{p^*_b}}\leq C\left(\int_{\mathbb{R}^N}\frac{1}{|x|^{ap}}|\nabla \varphi|^p\,dx\right)^{1/p}\qquad\hbox{ for all }\; \varphi\in C^\infty_c(\mathbb{R}^N),
$$
where $p^*_b=\frac{Np}{N-p(1+a-b)}$.
\end{lemma}

The Hardy   inequality is  a particular case of those established in \cite{CKN1984}.

\begin{proposition}[Hardy's inequality]
\label{HI}
Let $p\in (1,N)$. There {exists} a constant $\mathcal{H}>0$ such that
\begin{equation}\label{HIe}
\int_{\mathbb{R}^N}\frac{1}{|x|^p}|u|^p\,dx\leq \frac{1}{\mathcal{H}^p}\int_{\mathbb{R}^N}|\nabla u|^p\,dx\quad\hbox{ for all }\; u\in W^{1,p}(\mathbb{R}^N).
\end{equation}
Furthermore, the constant $\mathcal{H}=\frac{N-p}{p}$ is optimal.
\end{proposition}

The proof of this proposition can be found in \cite{AzoreroAlonso1998}. Observe that this inequality can be extended for $p\geq1$ by application of Lemma \ref{CKN}, for $b=1$ and $a=0$, by a density argument  $C^\infty_c(\mathbb{R}^N)$ is dense in $W^{1,p}(\mathbb{R}^N)$. 

The following version of Caffarelli-Kohn-Nirenberg in $BV(\Omega)$ is also established. 
\begin{proposition}
\label{PI}
Let $0\leq a< N-1 $ and $a\leq b\leq 1+a$. There exists a constant $C>0$ such that the following inequality is established
$$
\left(\int_{\mathbb{R}^N}\frac{1}{|x|^{b1^*_a}}|u|^{1^*_a}\,dx\right)^\frac{1}{1^*_a}\leq C\int_{\mathbb{R}^N}\frac{1}{|x|^a}|Du|\quad\hbox{ for all }\;u\in BV(\mathbb{R}^N),
$$
where $1^*_a=\frac{N}{N-(1+a-b)}$.
\end{proposition}

The proof of the above proposition can be found in \cite{OrtizHiroshi2024}.
\begin{remark}

Observe that Proposition \ref{PI}, for $a=b=0$, reduces to  Sobolev inequality, whose proof can be found in \cite[Theorem 3.47]{AmbrosioFuscoPallara}. Moreover, the constant
 \begin{equation}
\label{sn}
S_N=\frac{1}{NC_N^{1/N}}
\end{equation}
is the best in this case (see \cite{Talenti}).

\end{remark}

An immediate consequence of Proposition \ref{SL} (see \cite{pe,Alvino1977,ct,crt}) is the following.
\begin{proposition}
\label{LSI}
There exists a constant $\gamma_1>0$ such that
\begin{equation*}
\|u\|_{L^{1^*, 1}(\mathbb{R}^N)}\leq  \gamma_1 \int_{\mathbb{R}^N}|Du| \quad\mbox{ for all }\; u\in BV(\mathbb{R}^N).
\end{equation*}
The better constant is given by
\begin{equation}
\label{Tg}
 \gamma_1=\frac{1}{(N-1)C_N^{1/N}}.
\end{equation}
\end{proposition}

\smallskip

\subsection{ Some fine properties of $BV$ functions}

Recall that, for $u\in L^1(\Omega)$, $u$ has an approximate limit at $x\in \Omega$ if there exists $\tilde{u}(x)$ such that 
$$
\lim_{\rho\downarrow0} \dashint_{B_{\rho}(0)}|u(y)-\tilde{u}(x))|\,dx=0,
$$
where $\dashint_{E}u\,dx=\frac{1}{|E|}\int_{E}u\,dx$. The set where this property does not hold is denoted by $S_u$. This is a $\mathcal{L}^N-$negligible Borel set \cite[Proposition 3.64]{AmbrosioFuscoPallara} . We say that $x$ is an approximate jump point of $u$ if there exist $u^+(x)\neq u^{-}(x)$ and $ \nu\in S^{N-1}$  such that
\begin{eqnarray*}
\lim_{\rho \downarrow 0}\dashint_{B^+_\rho(x, \nu) }|u(y)-u^+(x)|\,dx&=&0\\
\lim_{\rho \downarrow 0}\dashint_{B^{-}_\rho(x, \nu) }|u(y)-u^{-}(x)|\,dx&=&0,
\end{eqnarray*}
where
\begin{eqnarray*}
B^+_\rho (x, \nu)&=&\{y\in B_\rho(x)\>:\> \langle y-x, \nu\rangle >0\}\\
B^{-} _\rho (x, \nu)&=&\{y\in B_\rho(x)\>:\> \langle y-x, \nu\rangle <0\}.
\end{eqnarray*}
The set of approximate jump points is denoted by $J_u$. The set $J_u$ is a Borel subset of $S_u$ \cite[Proposition 3.69]{AmbrosioFuscoPallara} and $\mathcal{H}^{N-1}(S_u \setminus J_u) = 0$.

For $u\in L^1(\Omega)$, $u^*\>:\>\Omega\setminus(S_u\setminus J_u)\to \mathbb{R}$ is called the canonical representative of u if 
\begin{equation}
\label{CanRep}
u^*(x)=
\left\{
\begin{array}{lcr}
\tilde{u}(x)&\mbox{ if }& x\in \Omega\setminus S_u\\
\frac{u^+(x)+u^{-}(x)}{2}& \mbox{ if }& x\in J_u.
\end{array}
\right.
\end{equation}

\begin{proposition}
\label{CRp}
Let $u \in BV(\Omega)$. The mollified functions $u*\rho_\epsilon$ pointwise converge to $u^*$ in its domain.
\end{proposition}

The proof of the above proposition as well as to more details about space $BV(\Omega)$ can be found in \cite{AmbrosioFuscoPallara} (to which we refer, in general, for further standard  notations not explicitly recalled here), see also \cite{AttouchButtazzoMichaille2006, EvansGariepy2015}.

\subsection{An  Anzellotti-Chen-Frid type theory of pairings}\label{anze}

In this section, we recall an extension of Anzellotti's theory \cite{Anzellotti1983} (see also \cite{CF}), given in \cite{Caselles2011}.

We set the space
$$
X_{\mathcal{M}}(\Omega)=\{{\bf z}\in L^{\infty}(\Omega,\mathbb{R}^N)\>:\> \mbox{div}\,{\bf z} \; \mbox{ is a finite Radon measure in }\;\Omega\}.
$$

We define
$$  
\langle \mbox{div}\,{\bf z},u\rangle=\int_{\Omega}u^*\,\mbox{div}\,{\bf z}\quad\hbox{ for all }\; {\bf z}\in X_{\mathcal{M}}(\Omega),\;\hbox{ for all }\; u\in BV(\Omega)\cap L^\infty(\Omega),
$$
where $u^*$ is as in \eqref{CanRep}.

For each $u\in W^{1,1}(\Omega)\cap L^\infty(\Omega)$ and ${\bf z}\in X_{\mathcal{M}}(\Omega)$ we define
\begin{equation*}
\langle {\bf z}, u\rangle_{\partial\Omega}=\int_{\Omega}u^*\,\mbox{div}{\bf z}+\int_{\Omega}{\bf z}\cdot \nabla u\,dx.
\end{equation*}

\begin{definition}
For each ${\bf z}\in X_{\mathcal{M}}(\Omega)$ and $u\in BV(\Omega)$, we define
$$
\langle {\bf z}, u\rangle_{\partial\Omega}=\langle {\bf z}, w\rangle_{\partial\Omega},
$$
where $w\in W^{1,1}(\Omega)\cap L^\infty(\Omega)$ is such that $ w = u $ on $ \partial\Omega $. 
\end{definition}

The above definition is well defined since the trace operator from $BV(\Omega)$ in $L^1(\partial\Omega)$ is onto and \cite[Lemma 5.5.]{Anzellotti1983} holds. On the other hand, arguing as in the proof of \cite[Theorem 1.2]{Anzellotti1983} we get the following result.

\begin{proposition}
\label{Prop2}
There exists a linear operator $ \beta\>:\>X_{\mathcal{M}}(\Omega)\to L^\infty(\partial\Omega)$ such that 
\begin{eqnarray*}
\label{E1}
\|\gamma_{\bf z}\|_{ L^\infty(\partial\Omega)}&\leq&\|{\bf z}\|_{ L^\infty(\Omega,\re^N)}\\
\langle {\bf z}, u\rangle_{\partial\Omega}&=&\int_{\partial\Omega}\beta_{\bf z}(x)u(x)\,d\mathcal{H}^{N-1}\quad \mbox{ for all }u\in BV(\Omega) \cap L^\infty(\Omega)\\
\label{E3}
\beta_{\bf z}(x)&=&{\bf z}(x)\cdot \nu(x)\quad\mbox{ for all }\; x\in \partial\Omega \quad\mbox{ if }\; {\bf z}\in C^1(\overline{\Omega},\mathbb{R}^N).
\end{eqnarray*}
\end{proposition}

The function $\beta_{\bf z}(x)$ is a weakly defined trace on $\partial\Omega$ of the normal component of ${\bf z}$, which will be denoted by $[{\bf z}, \nu]$. 
\begin{definition}
\label{zDu}
Let ${\bf z} \in X_{\mathcal{M}}(A)$ and $u\in BV(A)\cap L^\infty(A)$ for all open set $A\subset\subset\Omega$. We define a linear functional $({\bf z}, Du)\>:\>C^\infty_c(\Omega)\to \mathbb{R}$ as 
$$
\langle ({\bf z}, Du),\varphi\rangle=-\int_{\Omega}u^*\varphi\,\mbox{div}\,{\bf z}-\int_{\Omega}u {\bf z}\cdot \nabla \varphi\,dx.
$$
\end{definition}

\begin{proposition}
\label{Dist}
Let $A\subset \Omega$ be a open set, ${\bf z} \in X_{\mathcal{M}}(A)$,  and $u\in BV(A)\cap L^\infty(A)$. Then  for all functions $\varphi\in C_c(A)$ the following inequality holds                                                                                                                                                                                                                                                                                                                                                                                                                
$$
|\langle ({\bf z}, Du), \varphi\rangle |\leq  \max_{A}|\varphi|\|{\bf z}\|_{L^{\infty}(A,\re^N)}\int_{A}|Du\,|,
$$
which means that $({\bf z}, Du)$ is a Radon measure in $\Omega$. Furthermore, $({\bf z}, Du)$, $|({\bf z}, Du)|$ are absolutely continuous with respect to the measure $|Du|$ in $\Omega$ and 
$$
\left|\int_{B}({\bf z}, Du)\right|\leq \int_{B} |({\bf z}, Du)|\leq \|{\bf z}\|_{L^\infty(A,\re^N)}\int_{B}|Du|\,,
$$
for all Borel sets $B$ and for all opens sets $A$ such that $B\subset A\subset \Omega$.
\end{proposition}

\begin{proposition}
\label{Mea}
Let ${\bf z}\in X_{\mathcal{M}}(\Omega)$ and $u\in BV(\Omega)\cap L^\infty(\Omega)$. Then the following identity holds
\begin{equation}\label{pairingb}
\int_{\Omega}u^*\, \mbox{div}\,{\bf z}+\int_{\Omega}({\bf z}, Du)=\int_{\partial\Omega}[{\bf z}, \nu]\,u\,d\mathcal{H}^{N-1}.
\end{equation}
\end{proposition}
Proofs of the propositions \ref{Dist} and \ref{Mea} can be found in \cite[Section 5]{Caselles2011} (see also \cite{Anzellotti1983}).

\begin{definition}
Let ${\bf z}\in X_{\mathcal{M}}(\Omega)$ and $u\in BV(\Omega)$ be such that $u^*\in L^1(\Omega, |\hbox{div}\, {\bf z}|)$. The linear functional $({\bf z}, Du)\>:\>C^\infty_c(\Omega)\to \mathbb{R}$ is defined by
\begin{equation}\label{pairing}
\langle ({\bf z}, Du), \varphi\rangle=-\int_{\Omega}\varphi u^*\hbox{div}\,{\bf z}-\int_{\Omega} u {\bf z}\cdot \nabla \varphi\,dx\quad\hbox{ for all }\, \varphi\in C^\infty_c(\Omega).
\end{equation}
\end{definition}

\begin{proposition}
\label{Tku}
Suppose that ${\bf z}\in X_{\mathcal{M}}(\Omega)$ and $u\in BV(\Omega)$ is such that $u^*\in L^1(\Omega, |\hbox{div}\,{\bf z}|)$ and $u^+, u^-\in L^1(J_u, |\hbox{div}\,{\bf z}|)$. Then
$$
({\bf z}, T_k(u))\to ({\bf z}, Du)\quad\hbox{ in }\; \mathcal{D}'(\Omega)\quad\hbox{ as }\; k\to \infty.
$$
\end{proposition}

The proof of this proposition can be found in \cite{OrtizPetitta2024}.

\section{Main  results for $L^N$ data}\label{main}

For the sake of exposition, in this section we restrict ourself to the case of a datum $f$ in $L^{N}(\Omega)$, while the extension to more general data will be discussed in Section \ref{sec5}. 

 We consider  
\begin{equation}
\label{P} 
\left\{
\begin{array}{rclc}
\displaystyle-\Delta_1u&=&\lambda\frac{D u}{|D u|}\cdot \frac{x}{|x|^2} +f&\quad \mbox{ in } \Omega,\\
u&=&0 &\quad \mbox{ on }\partial\Omega,
\end{array}
\right.
\end{equation} 
where  $\Omega\subset \mathbb{R}^N$ $(N\geq 2)$ is an open set with bounded Lipschitz boundary containing the origin, $0<|\lambda|<N-1$, and $f\in L^N(\Omega)$ is a not identically null function that satisfies
\begin{equation}
\label{H}
{S_N}\|f\|_{L^N(\Omega)}+\frac{|\lambda|}{N-1}\leq 1,
\end{equation}
where $S_N$ is as in \eqref{sn}. Note that, when $\lambda=0$, this threshold  coincides with the one  considered in \cite{ct,mst}. 
 
Here is how  solutions of  \eqref{P} are intended as a natural extension  of the concept of solution introduced in \cite{acm2001, AndreuMazonMollCaselles2004, AndreuCasellesMazon2004}. 

 \begin{definition}
\label{DefSol}
We say that $u\in BV(\Omega)$ is a solution of \eqref{P} if  there exists ${\bf z}\in X_{\mathcal{M}}(\Omega)$ with $\|{\bf z}\|_{L^\infty(\Omega,\mathbb{R}^N)}\leq 1$ such that 
\begin{itemize}
\item[$(1)$] $\displaystyle-\mbox{\rm div}\,{\bf z}=\lambda{\bf z}\cdot\frac{x}{|x|^2}+f\quad \mbox{ in }\;\mathcal{D}'(\Omega)$;
\item[$(2)$] $({\bf z}, DT_k(u))=|DT_k(u)|\quad\mbox{ as Radon measures in }\;\Omega,\;\hbox{ for all }\,k>0$;
\item[$(3)$] $ [{\bf z}, \nu](x)\in \mbox{\rm Sgn}(-u(x))\quad \mathcal{H}^{N-1}-\mbox{ a.e. }
\;x\in \partial\Omega.$
\end{itemize}
\end{definition}

Our main result is the following.
\begin{theorem} \label{exi} Let $f\in L^N(\Omega)$. If $f$ satisfies \eqref{H}, then there exists a solution of problem \eqref{P} in the sense of Definition \ref{DefSol}.
\end{theorem}

As in \cite{ct,mst}, the key  argument in order to prove  Theorem \ref{P} consists in the  study of the asymptotic behaviour,  as $p\to 1^+$,  of  the solutions $u_p$ of the problem
\begin{equation}
\label{P1}
\left\{
\begin{array}{rclc}
\displaystyle-\Delta_p u&=&\lambda|\nabla u|^{p-2}\nabla u\cdot \frac{x}{|x|^2}+f&\mbox{ in  }  \Omega,\\
u&=&0&\mbox{ on }\partial\Omega,
\end{array}
\right.
\end{equation}
for which a solution does exist  by \cite[Theorem 3.1]{LeonoriPetitta2007}. 
More precisely, we shall find uniform  estimates with respect to  $p$ in order to pass the limit as $p\to 1^+$.
To do that  the main difficulty relies  in the accurate  analysis  of a priori estimates and convergence properties  of the family
 $$
\left(\lambda|\nabla u_p|^{p-2}\nabla u_p\cdot \frac{x}{|x|^2}\right)_{p>1},
$$
which is overcome, as in \cite{OrtizPetitta2024}, by mean  of Hardy inequality \eqref{HIe} as well as a suitable   truncation argument.

The asymptotic behaviour of the family $(u_p)_{p>1}$ provided  $f$ satisfies \eqref{H} is the content of the following result which is in sharp continuity with the previous results of \cite{ct, mst, dgop,OrtizPetitta2024}.

 \begin{theorem}
\label{AsympB}
Assume that $f$ satisfies 
$$
{S_N}\|f\|_{L^N(\Omega)}+\frac{|\lambda|}{N-1}< 1.
$$
 Then
$$
u_p\to 0\qquad \mbox{ a.e. in }\quad\Omega\quad\mbox{ as }\quad p\to 1^+.
$$
If  $f$ satisfies 
$$
{S_N}\|f\|_{L^N(\Omega)}+\frac{|\lambda|}{N-1}= 1,
$$
then  there exist $u\in BV(\Omega)$ such that 
$$
u_p\to u\qquad\mbox{ in }\quad L^s(\Omega),\; s\in [1,1^*)\quad\mbox{ as } \quad p\to 1^+.
$$
\end{theorem}

Furthermore, we shall show the following  $L^\infty(\Omega)$ regularity result for  solutions to problem \eqref{P}  whose proof is a re-adaptation of the classical Stampacchia's argument to our framework ({see} also \cite[Theorem 3.5]{MS2013}). 

\begin{theorem}
\label{bo}
Let $f\in L^q(\Omega)$ with $q>N$. Assume that $f$ satisfies \eqref{H}. Then the solutions $u_p$ of \eqref{P1} satisfies 
\begin{equation}
\|u_p\|_{L^{\infty}(\Omega)}\leq C, \label{linf}
\end{equation}
for some constant $C>0$ not depending on $p$. 
In particular the limit  of $u_p$, that is   $u$ solution of  \eqref{P} belongs to $L^\infty(\Omega)$. 
\end{theorem}

\subsection{Some preliminary  a priori estimates}
\label{F} 
Assume that there exists $\bar p\in (1, N)$ such that 
\begin{equation}
\label{EstLbd}
|\lambda|<\frac{N-p}{p}\qquad \hbox{ for all }\; p\in (1,\bar p].
\end{equation}

Observe that, once \eqref{H} is in force, for $f\not\equiv 0$, then \eqref{EstLbd} always hold for some $\bar p$.

Now, fixed $f\in L^N(\Omega)$, consider for any  $p\in(1, \bar p]$ the problem
\begin{equation*}
\label{Lnp}
\left\{
\begin{array}{rclr}
-\Delta_p u &=&\lambda|\nabla u|^{p-2}\nabla u\cdot\frac{x}{|x|^2}+f&\mbox{ in  }  \Omega,\\
u&=&0&\mbox{ on }\partial\Omega.
\end{array}
\right.
\end{equation*}
By \cite[Theorem 3.1 {(i)}]{LeonoriPetitta2007}, there exists a weak solution $u_p\in W^{1,p}_0(\Omega)\cap L^\infty(\Omega)$, i.e., $u_p$ satisfies
\begin{equation}
\label{WSol}
\int_{\Omega}|\nabla u_p|^{p-2}\nabla u_p\cdot\nabla \varphi\,dx=\lambda\int_{\Omega}|\nabla u_p|^{p-2}\nabla u_p\cdot \frac{x}{|x|^2}\varphi\,dx+\int_{\Omega}f\varphi\,dx\quad\forall\,\varphi\in W^{1,p}_0(\Omega).
\end{equation}
Observe that the first term on the right-hand side of \eqref{WSol} is well defined because of the  Hardy inequality \eqref{HIe}; in fact, using the H\"older inequality first and then \eqref{HIe}, one has 
$$
\left|\int_{\Omega}\frac{|\nabla u|^{p-2}\nabla u\cdot x}{|x|^2}\varphi\,dx\right|\leq \left(\int_{\Omega} |\nabla u|^p dx \right)^{1/p'}\left(\int_{\Omega} \frac{|\varphi|^p}{|x|^p}\right)^{1/p}\leq \frac{p}{N-p} \|u\|_{W^{1,p}_0(\Omega)}^{p-1} \|\varphi\|_{W^{1,p}_0(\Omega)}.
$$
Moreover, using $u_p$ as test function in \eqref{WSol}, and using \eqref{HIe} as above,  we get

$$
\left(1-\frac{|\lambda| p}{N-p}\right)\int_{\Omega} |\nabla u_p|^{p} dx\leq \|f\|_{L^{N}(\Omega)}{\|u\|_{L^{1^*}(\Omega)}}
$$
and so, using Sobolev and {H\"older inequalities}, one has 
\begin{equation} 
\label{B}
\|u_p\|_{W^{1,p}_0}\leq \left[\frac{\|f\|_{L^N(\Omega)}S_N}{1-\frac{|\lambda| p }{N-p}}\right]^{1/(p-1)}|\Omega|^{1/p}\qquad\mbox{ for all }\; p\in (1, \bar p].
\end{equation}

 The following result is a consequence of  the above inequality and condition \eqref{H}.
 
\begin{lemma}
\label{UB}
 Let $f\in L^{N}(\Omega)$ and $\lambda\in \re$ be such that   \eqref{H} holds.  Then there exists $u \in BV (\Omega)$ such that, up  to subsequences,
 \begin{eqnarray*}
 \label{ConLp}
 u_p\to u&\mbox{in }\quad L^s(\Omega),\; s\in [1, 1^*)\\
 \nonumber u_p\rightharpoonup u&\mbox{ in } \quad L^{1^*}(\Omega)\\
 \nonumber \nabla u_p\rightharpoonup Du&*-\hbox{weakly in}\quad\mathcal{M}(\Omega,\re^N),
 \end{eqnarray*}
 as $ p\to 1^+$.
\end{lemma}

\dem
Note that
$$
\lim_{p\to 1^+} \left[\frac{1-\frac{|\lambda|}{N-1}}{1-\frac{|\lambda| p}{N-p}}\right]^{1/(p-1)}=e^{|\lambda| N/[(N-1)(N-1-|\lambda|)]}.
$$
Hence  using \eqref{H}  in \eqref{B},  we deduce that there exists $C>0$ such that
$$
\int_{\Omega}|\nabla u_p|^p\,dx\leq C|\Omega|\quad\hbox{ for all }\;p\in (1, \bar p].
$$
By Young's inequality, we have that
$$
\int_{\Omega}|\nabla u_p|\,dx\leq \frac{1}{p}\int_{\Omega}|\nabla u_p|^p\,dx+\frac{p-1}{p}|\Omega|\leq(C+1) |\Omega|\quad\hbox{ for all }\; p\in (1, \bar p].
$$
Since $u_p=0$ on $\partial\Omega$, the family $(u_p)_{p>1}$ is uniformly bounded in $BV(\Omega)$. Consequently, by  compact embedding in $BV$ (see for instance  \cite[Proposition 3.13]{AmbrosioFuscoPallara}), there exists $u\in BV(\Omega)$ such that, up to subsequences,
\begin{equation*}
\label{ConvUp}
u_p\to u\quad\mbox{in}\quad L^s(\Omega),\; s\in [1, 1^*)
\end{equation*}
and 
$$
\nabla u_p\rightharpoonup Du\qquad*-\hbox{weakly in }\mathcal{M}(\Omega,\re^N), 
$$
as $p\to 1^+$. Recall that, by  the Sobolev inequality, one also gets 
 $$
u_p\rightharpoonup u\quad\mbox{ in }\;L^{1^*}(\Omega)\quad\mbox{ as }\; p\to 1^+.
$$
\fim

\begin{lemma}
\label{ConvZ}
Let $q\geq 1$. There exists ${\bf z}\in L^\infty(\Omega,\mathbb{R}^N)$ with $\|{\bf z}\|_{L^\infty(\Omega,\mathbb{R}^N)}\leq 1$ such that 
$$
|\nabla u_p|^{p-2}\nabla u_p\rightharpoonup{\bf z}\quad \mbox{ in }\;L^q(\Omega,\mathbb{R}^N)\quad\mbox{ as }\;p\to 1^+.
$$
\end{lemma}			

The proof of the previous lemma is standard and can be found, for instance,  in \cite[Lemma 4.3]{OrtizPetitta2024}. 

\section{Proof of  main theorems}\label{mains}

In this section, we shall prove  Theorem \ref{exi}, Theorem \ref{AsympB}, and Theorem \ref{bo}.

\subsection{Proof of Theorem \ref{exi}}
We  show that the limits $u$ and ${\bf z}$ detected in Lemma \ref{UB} and \ref{ConvZ}, respectively, satisfy $(1)-(3)$ of Definition \ref{DefSol}.

{\it Proof of $(1)$.} In order to prove $(1)$  we  pass to the limit in \eqref{WSol}.

Let $\varphi\in C^\infty_c(\Omega)$. By Lemma \ref{ConvZ}, it follows that
$$
\lim_{p\to1^+}\int_{\Omega}|\nabla u_p|^{p-2}\nabla u_p\cdot\nabla \varphi\,dx=\int_{\Omega}{\bf z}\cdot \varphi\,dx
$$
and
$$
\lim_{p\to 1^+}\int_{\Omega}|\nabla u_p|^{p-2}\nabla u_p\cdot\frac{x}{|x|^2}\varphi\,dx=\int_{\Omega}{\bf z}\cdot\frac{x}{|x|^2}\varphi\,dx.
$$

Letting $p\to1^+$ in \eqref{WSol}, we obtain that
 \begin{equation*}
 \label{ConvS1}
 \int_{\Omega}{\bf z}\cdot\nabla \varphi\,dx=\lambda\int_{\Omega}{\bf z}\cdot\frac{x}{|x|^2}\varphi\,dx+\int_{\Omega}f\varphi\,dx,
 \end{equation*}
 i.e. 
 \begin{equation}\label{(1)}
 \displaystyle-\mbox{\rm div}\,{\bf z}=\lambda{\bf z}\cdot\frac{x}{|x|^2}+f\quad \mbox{ in }\;\mathcal{D}'(\Omega).
 \end{equation}
Notice, in particular, that  $\hbox{div}\,{\bf z}$ is not only a measure in $\mathcal{M}(\Omega)$; in fact,   by \eqref{(1)} one has that 
  $
  \hbox{div}\,{\bf z} \ \ \in L^q(\Omega), \ \  \forall\, q\in [1, N),  
  $
or, more precisely
$$
 \hbox{div}\,{\bf z} \ \  \in L^{N,\infty}(\Omega).
$$

{\it Proof of  $(2)$.} Note that, by {Proposition \ref{Dist}}, we have that
$$
({\bf z}, DT_k(u))\leq |DT_k(u)|\quad\hbox{ as Radon measures on }\;\Omega.
$$

In order to show the reverse inequality, let $\varphi\in C^\infty_c(\Omega)$ be a nonnegative function, we use $\varphi T_k(u_p)$ as test function in \eqref{WSol} to obtain
\begin{align*}
&\int_{\Omega}\varphi |\nabla T_k(u_p)|^p\,dx+\int_{\Omega}T_k(u_p)|\nabla u_p|^{p-2}\nabla u_p\cdot \nabla \varphi\,dx\\
&\qquad=\lambda\int_{\Omega}|\nabla u_p|^{p-2}\nabla u_p\cdot\frac{x}{|x|^2}\varphi T_k(u_p)\,dx+\int_{\Omega}f\varphi T_k(u_p)\,dx.
\end{align*}
Using Young's inequality, we have that
\begin{align}
\int_{\Omega}\varphi |\nabla T_k(u_p)|&\leq\frac{1}{p}\int_{\Omega}\varphi|\nabla T_k(u_p)|^p\,dx+\frac{p-1}{p}\int_{\Omega}\varphi\,dx\nonumber\\
&\leq \frac{\lambda }{p}\int_{\Omega}|\nabla u_p|^{p-2}\nabla u_p\cdot \frac{x}{|x|^2}\varphi T_k(u_p)\,dx+\frac{1}{p}\int_{\Omega}f\varphi T_k(u_p)\,dx \label{E} \\
\nonumber&-\frac{1}{p}\int_{\Omega}T_k(u_p)|\nabla u_p|^{p-2}\nabla u_p\cdot\nabla \varphi\,dx+\frac{p-1}{p}\int_{\Omega}\varphi\,dx.
\end{align}
Notice that, as $\displaystyle \varphi T_k(u_p)x/|x|^2\to \varphi T_k(u)x/|x|^2$ in $L^q(\Omega, \mathbb{R}^N)$ and $|\nabla u_p|^{p-2}\nabla u_p\rightharpoonup {\bf z}$ in $L^{q'}(\Omega,\mathbb{R}^N)$, for $q\in (1, N)$, it follows that
\begin{equation}
\label{Ck}
\lim_{p\to 1^+}\int_{\Omega}|\nabla u_p|^{p-2}\nabla u_p\cdot \frac{x}{|x|^2}\varphi T_k(u_p)\,dx=\int_{\Omega}{\bf z}\cdot \frac{x}{|x|^2}\varphi T_k(u)\,dx.
\end{equation}
Similarly, one can prove that
$$
\lim_{p\to 1^+}\int_{\Omega}T_k(u_p)|\nabla u_p|^{p-2}\nabla u_p\cdot \nabla \varphi\,dx=\int_{\Omega} T_k(u) {\bf z}\cdot \nabla \varphi\,dx.
$$
Consequently, passing to limit in \eqref{E} as $p\to 1^+$, since $v\mapsto \int_{\Omega}\varphi |Dv|$ is lower semicontinuous  with respect to the $L^1$ convergence, 
\begin{align}
\label{EsTku}
\int_{\Omega}\varphi |D T_k(u)|&\leq \lambda\int_{\Omega} {\bf z}\cdot\frac{x}{|x|^2}\varphi T_k(u)\,dx+\int_{\Omega}f\varphi T_k(u)\,dx-\int_{\Omega}T_k(u) {\bf z}\cdot \nabla \varphi\,dx\nonumber\\ 
&\stackrel{\eqref{(1)}}{=} -\int_{\Omega}\varphi T_k(u)\hbox{div}\,{\bf z}\,dx-\int_{\Omega}T_k(u){\bf z}\cdot \nabla \varphi\,dx\\
\nonumber&\stackrel{\eqref{pairing}}{=} \int_{\Omega} ({\bf z}, DT_k(u)) \varphi\,,
\end{align}
for all $0\leq\varphi\in C^\infty_c(\Omega)$, and so, the reverse inequality holds. 

{\it Proof of $(3)$.} Notice that, by Proposition \ref{Prop2},
\begin{equation}
\label{Nn}
T_k(u)[{\bf z}, \nu]+|T_k(u)|\geq 0\quad\mathcal{H}^{N-1}-\hbox{ a.e. in }\;\partial\Omega.
\end{equation}
In order to prove the reverse inequality, use $T_k(u_p)$ as test function in \eqref{WSol} to obtain
$$
\int_{\Omega}|\nabla T_k(u_p)|^p\,dx=\lambda\int_{\Omega}|\nabla u_p|^{p-2}\nabla u_p\cdot\frac{x}{|x|^2}T_k(u_p)\,dx+\int_{\Omega}f T_k(u_p)\,dx.
$$
As $u_p=0$ on $\partial\Omega$, using Young inequality, we have that
\begin{align}
\label{kp}
\nonumber\int_{\Omega}|\nabla T_k(u_p)|&+\int_{\partial\Omega}|T_k(u_p)|\,d\mathcal{H}^{N-1}\\
&\leq \frac{1}{p}\int_{\Omega}|\nabla T_k(u_p)|^p\,dx+\frac{p-1}{p}|\Omega|\\
\nonumber&=\frac{\lambda}{p}\int_{\Omega}|\nabla u_p|^{p-2}\nabla u_p\cdot\frac{x}{|x|^2}T_k(u_p)\,dx+\frac{1}{p}\int_{\Omega} fT_k(u_p)\,dx+\frac{p-1}{p}|\Omega|.
\end{align}
Note that, as in \eqref{Ck}, one also can prove that
$$
\lim_{p\to 1^+}\int_{\Omega}|\nabla u_p|^{p-2}\nabla u_p\cdot \frac{x}{|x|^2}T_k(u_p)\,dx=\int_{\Omega}{\bf z}\cdot\frac{ x}{|x|^2}T_k(u)\,dx.
$$
Letting $p\to 1^+$ in \eqref{kp}, using the lower semi-continuity, we obtain that
\begin{align*}
\label{EsTk}
\int_{\Omega}|DT_k(u)|+\int_{\partial\Omega}|T_k(u)|\,d\mathcal{H}^{N-1}&\leq \lambda\int_{\Omega}{\bf z}\cdot \frac{x}{|x|^2} T_k(u)\,dx+\int_{\Omega} f T_k(u)\,dx\\
&\stackrel{\eqref{(1)}}{=}-\int_{\Omega}T_k(u)\hbox{div}\,{\bf z}\,dx\\
&\stackrel{\eqref{pairingb}}{=}\int_{\Omega}({\bf z}, DT_k(u))-\int_{\partial\Omega}[{\bf z},\nu]T_k(u)\,d\mathcal{H}^{N-1}.
\end{align*}
Since $({\bf z}, DT_k(u))=|DT_k(u)|$,
\begin{equation}
\label{Np}
\int_{\partial\Omega}|T_k(u)|\,d\mathcal{H}^{N-1}+\int_{\partial\Omega}[{\bf z}, \nu] T_k(u)\,d\mathcal{H}^{N-1}\leq 0.
\end{equation}

Thus, from \eqref{Nn} and \eqref{Np},
\begin{equation}
\label{Idk}
\int_{\partial\Omega}(|T_k(u)|+T_k(u)[{\bf z}, \nu])\,d\mathcal{H}^{N-1}=0.
\end{equation}
Since $T_k(u)\to u$ in $L^1(\Omega)$ and $\displaystyle \int_{\Omega}|DT_k(u)|\leq \int_{\Omega}|Du|$ for all $k\geq 1$, $(T_k(u))_{k\geq1}$ strictly converges to $u$. But this, in turn, implies, by continuity of trace operator with respect to the strict convergence,  that
$$
T_k(u)\to u\quad\hbox{ in }\; L^1(\partial\Omega),\quad\mathcal{H}^{N-1}-\hbox{ a.e. in  }\;\partial\Omega\quad\hbox{ as }\; k\to \infty.
$$

Letting $k\to \infty$ in \eqref{Idk}, we obtain that
$$
\int_{\partial\Omega}(|u|+u[{\bf z}, \nu])\,d\mathcal{H}^{N-1}=0.
$$

Therefore $u$ and ${\bf z}$ as in Lemma \ref{UB} and Lemma \ref{ConvZ} respectively, satisfy $(1)-(3)$ of Definition \ref{DefSol}.

\fim\bk
\begin{remark}\label{remk}
If we assume that the limit $u$ in Lemma \ref{UB} is such that $ u^*\in L^1(\Omega, |\hbox{div}\,{\bf z}|)$ and $u^+, u^-\in L^1(\Omega, J_u)$, then passing to the limit in \eqref{EsTku}, by Proposition \ref{Tku}, and using the weak lower semicontinuity, 
$$
\int_{\Omega}\varphi|Du|\leq \int_{\Omega} \varphi({\bf z}, Du) \quad\forall\; 0\leq \varphi\in C^\infty_c(\Omega).
$$
On the other hand, by Proposition \ref{Dist} and Proposition \ref{Tku}, we have that
$$
\int_{\Omega}\varphi({\bf z}, Du)=\lim_{k\to\infty}\int_{\Omega}\varphi({\bf z}, DT_k(u))\leq\limsup_{k\to\infty}\int_{\Omega}\varphi |DT_k(u)|=\int_{\Omega}\varphi|Du|\quad\forall\; 0\leq \varphi\in C^\infty_c(\Omega).
$$
Consequently,
$$
\int_{\Omega}\varphi ({\bf z}, Du)=\int_{\Omega}\varphi|Du|\quad\forall\, 0\leq \varphi\in C^\infty_c(\Omega).
$$

\end{remark}

 \subsection{Proof of Theorem \ref{AsympB}}
 
 We split the proof of Theorem  \ref{AsympB} into  two cases:
 
 {\bf Case 1.} Suppose that
 $$
 S_N\|f\|_{L^N(\Omega)}+ \frac{|\lambda|}{N-1} <1.
 $$
 Notice that, given $\epsilon>0$,  there exist $\hat p>1$ such that 
 $$
 \frac{ \|f\|_{L^N(\Omega)}  S_N}{1-\frac{|\lambda| p}{N-p}}<1-\epsilon\quad \hbox{ for all }\; p\in (1, \hat p].
 $$
 Thus, from \eqref{B}, 
$$
u_p\to 0\quad\hbox{ in } \; W^{1,p}_0(\Omega)\quad\hbox{ as } p\to 1^+
$$
and so, up  to subsequences,
$$
u_p(x)\to 0\quad\hbox{ a.e. }\; x\in \Omega\quad\hbox{ as }\; p\to 1^+.
$$
 
{\bf Case 2.} Assume that
$$
 S_N\|f\|_{L^N(\Omega)}+\frac{|\lambda|}{N-1}=1.
 $$
Notice that
$$
\lim_{p\to 1^+}\left(\frac{\|f\|_{L^N(\Omega)} S_N }{1-\frac{|\lambda| p}{N-p}}\right)^\frac{p}{p-1}=\lim_{p\to1^+}\left(\frac{1-\frac{|\lambda|}{N-1}}{1-\frac{|\lambda| p}{N-p}}\right)^{p/(p-1)}=e^{|\lambda| N/[(N-1)(N-1-|\lambda|)]}.
$$
But this, in \eqref{B}, implies that there exists $C>0$ such that
$$
\int_{\Omega}|\nabla u_p|^p\,dx\leq C|\Omega|\quad \mbox{ for all }\; p\in (1, \bar p].
$$
Using  the Young inequality, we have that 
\begin{equation}
\int_{\Omega}|\nabla u_p|\,dx\leq \frac{1}{p}\int_{\Omega}|\nabla u_p|^p\,dx+\frac{p-1}{p}|\Omega|\leq (C+1)|\Omega|\quad \mbox{ for all }\; p\in (1, \bar p].\label{bve}
\end{equation}
Hence, as $u_p=0$ on $\partial\Omega$, $(u_p)_{p>1}$ is uniformly bounded in $BV(\Omega)$ and so, there exists $u\in BV(\Omega)$ such that, up  to subsequences, 
$$
u_p\to u\quad\hbox{ in }\; L^s(\Omega)\;\forall \,s\in [1, 1^*)\quad\hbox{ as }\;p\to 1^+.
$$
This complete the proof.
\fim

\subsection{Proof of Theorem \ref{bo}.}  
We use standard Stampacchia's type argument (\cite{Stampacchia1965}). We define
$$
A_{k,p}=\{x\in \Omega \colon |u_p(x)|\geq k\}.
$$
Use $G_k(u_p)$ as test function in \eqref{WSol} to obtain
$$
\int_{\Omega}|\nabla G_k(u_p)|^p\,dx=\lambda\int_{\Omega}|\nabla G_k(u_p)|^{p-2}\nabla G_k(u_p)\cdot \frac{x}{|x|^2}G_k(u_p)\,dx+\int_{\Omega}fG_k(u_p)\,dx.
$$
Using H\"older's, Hardy's and Sobolev's inequalities, the latter  becomes
$$
\left(1-\frac{|\lambda| p}{N-p}\right)\int_{\Omega}|\nabla G_k(u_p)|^p\,dx\leq \|f\|_{L^q(\Omega)}S_N\left(\int_{\Omega}|\nabla G_k(u_p)|^p\,dx\right)^{1/p}|A_{k,p}|^{1/N-1/q+(p-1)/p}
$$ 
 and so
 \begin{equation}
 \label{EA}
 \int_{\Omega}|\nabla G_k(u_p)|^p\,dx\leq \left(\frac{\|f\|_{L^q(\Omega)}S_N}{1-\frac{|\lambda| p}{N-p}}|A_{k,p}|^{1/N-1/q}\right)^{p/(p-1)}|A_{k,p}|\quad\hbox{ for all }\; p\in (1, \bar p].
 \end{equation}
 By the  Young inequality and \eqref{EA}, we obtain that 
 $$
 \int_{A_{k,p}}|\nabla G_k(u_p)|\,dx\leq \frac{1}{p}\left(\frac{\|f\|_{L^q(\Omega)}S_N}{1-\frac{|\lambda| p}{N-p}}|A_{k,p}|^{1/N-1/q}\right)^{p/(p-1)}|A_{k,p}|+\frac{p-1}{p}|A_{k,p}|,
 $$
 for all $p\in (1, \bar p]$. Using Chebyshev's inequality one has 
 $$
 |A_{k,p}|\leq \frac{1}{k^{N/(N-1)}}\int_{\Omega}|u_p|^\frac{N}{N-1}\,dx\stackrel{\eqref{bve}}{\leq} \frac{1}{k^{N/(N-1)}}(S_N(C+1)|\Omega|)^{N/(N-1)}\quad\hbox{ for all }\; p\in (1, \bar p].
 $$
 Hence, there exists $k_0\geq 1$ such that
 $$
 \frac{\|f\|_{L^q(\Omega)}S_N}{1-\frac{|\lambda| p}{N-p}}|A_{k,p}|^{1/N-1/q}<1\quad\hbox{ for all } \; k\geq k_0\;\mbox{ for all }\; p\in (1, \bar p].
 $$
 
 Thus
 \begin{equation*}
 \int_{A_{k,p}}|\nabla G_k(u_p)|\,dx\leq |A_{k,p}|\quad\hbox{ for all }\; k\geq k_0,\;\mbox{ for all }\; p\in (1, \bar p].
 \end{equation*}
But this, in turn, implies that
 $$
\int_{A_{k,p}}|G_k(u_p)|\,dx\leq S_N |A_{k,p}|^{1+1/N}\quad\hbox{ for all }\; k\geq k_0,\;\mbox{ for all }\; p\in (1, \bar p]
 $$
and so, for $l\geq k\geq k_0$, 
 $$
 |A_{l,p}|\leq \frac{S_N}{l-k}|A_{k,p}|^{1+1/N}\quad\mbox{ for all }\; p\in (1, \bar p].
 $$
 
Consequently, by \cite[Lemma 4.1]{Stampacchia1965},
 $$
 |A_{k_0+d,p}|=0\quad\mbox{ for all }\; p\in (1, \bar p],
 $$
 where $d=2^{N+1}S_N|A_{k_0,p}|^{1/N}$, and so \eqref{linf}  holds; moreover, since $u_p\to u$ a.e. in $\Omega$ then  $u \in L^\infty(\Omega)$. 
\fim

\begin{remark}
An easy consequence of Theorem \ref{bo} (see also Remark \ref{remk}) is that,   if $f\in L^q(\Omega)$ with $q>N$ in Theorem \ref{exi}, then the solution $u$ satisfies 
\begin{itemize}
\item [$(2)'$] $({\bf z}, Du)=|Du|\quad\hbox{ as Radon measures on }\;\Omega$,
\end{itemize}
instead of $(2)$ in Definition \ref{DefSol}. 
\end{remark}

\section{Data $f$ in  $L^{N, \infty}(\Omega)$ and $W^{-1,\infty}(\Omega)$} \label{sec5}
In this section we show how the previous argument can be use to extend some of the results to the case of more general data.  
Let $f\in L^{N,\infty}(\Omega)$. Assume that $f$ satisfies 
\begin{equation}
\label{HL}
{ \gamma_1}\|f\|_{L^{N,\infty}(\Omega)}+\frac{|\lambda|}{N-1} \leq 1,
\end{equation}
where $\gamma_1$ is as in \eqref{Tg}.
Let $\bar p>1$. For each $p\in (1, \bar p]$, consider the problem 
\begin{equation}
\label{PL}
\left\{
\begin{array}{rclc}
-\Delta_p u&=&\lambda |\nabla u|^{p-2}\nabla u\cdot\frac{x}{|x|^2}+f&\hbox{ in }\;\Omega\\
   		u&=&0&\hbox{ on }\;\partial\Omega.
\end{array}
\right.
\end{equation}
Again (\cite{LeonoriPetitta2007})   problem \eqref{PL} has solution, say $u_p$, that is, $u_p$ satisfies 
\begin{equation}
\label{wsL}
\int_{\Omega}|\nabla u_p|^{p-2}\nabla u_p\cdot \nabla \varphi\,dx=\lambda\int_{\Omega}|\nabla u_p|^{p-2}\nabla u_p\cdot \frac{x}{|x|^2}\varphi\,dx+\int_{\Omega} f\varphi\,dx\quad\hbox{ for all }\;\varphi\in W^{1,p}_0(\Omega).
\end{equation}

Use $u_p$ as test function in \eqref{PL}, by H\"older and Hardy inequalities,   and also using Proposition \ref{LSI}, one gets 
\begin{align*}
\int_{\Omega}|\nabla u_p|^p\,dx\leq \frac{|\lambda| p}{N-p}\int_{\Omega}|\nabla u_p|^p\,dx+\gamma_1\|f\|_{L^{N,\infty}(\Omega)}\left(\int_{\Omega}|\nabla u_p|^p\,dx\right)^{1/p}|\Omega|^{(p-1)/p}
\end{align*}
or, by rearranging
$$
\int_{\Omega}|\nabla u_p|^p\,dx\leq \left(\frac{{\ \gamma_1} \|f\|_{L^{N,\infty}(\Omega)}}{1-\frac{|\lambda| p }{N-p}}\right)^{p/(p-1)}|\Omega|\quad\hbox{ for all }\; p\in (1, \bar p].
$$
Since, using \eqref{HL},
$$
\limsup_{p\to 1^+}\left(\frac{{ \gamma_1} \|f\|_{L^{N,\infty}(\Omega)}}{1-\frac{|\lambda| p}{N-p}}\right)^{p/(p-1)}<\infty,
$$
there exists $u\in BV(\Omega)$ such that, {up to subsequences},
 \begin{eqnarray*}
  \nabla u_p\rightharpoonup Du&*-\hbox{weakly in}\; \mathcal{M}(\Omega,\re^N)  \\
 u_p\to u&\hbox{ in }\; L^q(\Omega),\; \forall \, q\in [1,1^*),\quad\hbox{ a.e. in }\;\Omega\\
 u_p\rightharpoonup u&\hbox{ in } L^{1^*}(\Omega), 
 \end{eqnarray*}
as $p\to 1$. 

On the other hand, using Lemma \ref{ConvZ},  there exists ${\bf z}\in L^\infty(\Omega, \mathbb{R}^N)$ with $\|{\bf z}\|_{L^\infty(\Omega,\mathbb{R}^N)}\leq 1$ such that 
\begin{equation*}
\label{Cz}
|\nabla u_p|^{p-2}\nabla u_p\rightharpoonup {\bf z}\quad\hbox{ in }\;L^q(\Omega, \mathbb{R}^N)\quad\hbox{ for all }\;q\geq 1\quad\hbox{ as }\; p\to 1^+.
\end{equation*}
{
Letting $p\to 1^+$ in \eqref{wsL}, by the very same  argument of the  proof of Theorem \ref{exi}, we obtain that
$$
\int_{\Omega}{\bf z}\cdot \nabla \varphi\,dx=\lambda\int_{\Omega}{\bf z}\cdot \frac{x}{|x|^2}\varphi\,dx+\int_{\Omega}f\varphi\,dx\quad\hbox{ for all }\; \varphi\in C^\infty_c(\Omega),
$$
i.e.,
$$
-\hbox{div}\,{\bf z}={\bf z}\cdot \frac{x}{|x|^2}+f \quad\hbox{ in  }\;\mathcal{D}'(\Omega).
$$

 It is worth to recall that $-\hbox{div}\,{\bf z}$ belongs to $L^{N, \infty}(\Omega)\subset \mathcal{M}(\Omega)$ and so, the Anzellotti's theory is still available. 
Considering the previous observations, following the same technique of the proof of Theorem \ref{exi}, one can show that $u$ and ${\bf z}$ satisfy $(1)-(3)$ of Definition \ref{DefSol}. One can also prove the  analogous of  Theorem \ref{AsympB}  by systematically using Proposition \ref{HoI} instead of the classical  H\"older inequality and Propositions \ref{SL} and \ref{LSI} in place of  the classical  Sobolev's inequalities.

Summarizing, one has the following results.
\begin{theorem} \label{exiinf} Let $f\in L^{N,\infty}(\Omega)$. If $f$ satisfies \eqref{HL}, then there exists a solution of problem \eqref{P} in the sense of Definition \ref{DefSol}.
\end{theorem}

 \begin{theorem}
Assume that $f$ satisfies 
$$
{\gamma_1}\|f\|_{L^{N,\infty}(\Omega)}+\frac{|\lambda|}{N-1}< 1.
$$
 Then
$$
u_p\to 0\qquad \mbox{ a.e. in }\quad\Omega\quad\mbox{ as }\quad p\to 1^+.
$$
If  $f$ satisfies 
$$
{\gamma_1}\|f\|_{L^{N,\infty}(\Omega)}+\frac{|\lambda|}{N-1}= 1,
$$
then  there exist $u\in BV(\Omega)$ such that 
$$
u_p\to u\qquad\mbox{ in }\quad L^s(\Omega),\; s\in [1,1^*)\quad\mbox{ as } \quad p\to 1^+.
$$
\end{theorem}

\begin{remark}
Similarly,  by re-adapting the previous idea, when $f\in W^{-1,\infty}(\Omega)$, assuming that  
\begin{equation*}
\label{HDinf}
\frac{|\lambda|}{N-1}+\|f\|_{W^{-1,\infty}(\Omega)}\leq 1,
\end{equation*}
  results analogous to  Theorem \ref{exi} and  Theorem  \ref{AsympB} can be established straightforwardly.
\end{remark}

\section{A generic drift  in $L^{N,\infty}(\Omega)$}\label{sec6}
Consider the following problem
\begin{equation}
\label{P1F}
\left\{
\begin{array}{rclc}
-\Delta_1 u &=& \frac{Du}{|Du|}\cdot {\bf F}+f& \hbox{ in }\; \Omega\\
               u&=&0&\hbox{ on }\;\partial\Omega,
\end{array}
\right.
\end{equation}
where $\Omega\subset \mathbb{R}^N$ is an open set with bounded Lipschitz boundary, $f$, for simplicity,  belongs to $L^N(\Omega)$ (being the case   $L^{N,\infty}(\Omega)$ easily treatable as in Section \ref{sec5}, see Remark \ref{remextend} below),  and the measurable vector function ${\bf F}\colon \Omega\to \mathbb{R}^N$ is such that $|{\bf F}|$ belong to $L^{N,\infty}(\Omega)$. Moreover, we assume
\begin{equation}
\label{HFf}
\gamma_1\||{\bf F}|\|_{L^{N,\infty}(\Omega)}+S_N\|f\|_{L^N(\Omega)}\leq 1,
\end{equation}
where $S_N$ and $\gamma_1$ are as in \eqref{sn} and \eqref{Tg} respectively. Note  that, if \eqref{HFf} is in force and $f\neq 0$, then $\||{\bf F}|\|_{L^{N, \infty}(\Omega)}<1/\gamma_1$. 

Also observe  that, if ${\bf F}(x)=\lambda x/|x|^2$ and $|\lambda|<N-1$, then \eqref{HFf} reduces to  \eqref{H}.

\subsection{Approximating problems and basic a priori estimates}
First of all, note that, as $f\neq 0$,  by \eqref{HFf}, there exists $\bar p>1$ such that
$$
1-\gamma_p\||{\bf F}|\|_{L^{N, \infty}(\Omega)}>0\quad\hbox{ for all }\; p\in (1, \bar p].
$$
For each $p\in (1, \bar p]$, we  consider the following problem 
\begin{equation}
\label{PpF}
\left\{
\begin{array}{rclc}
-\Delta_p u &=& |\nabla u|^{p-2}\nabla u\cdot {\bf F}+f& \hbox{ in }\; \Omega\\
               u&=&0&\hbox{ on }\;\partial\Omega.
\end{array}
\right.
\end{equation}
 Reasoning as in  \cite[Theorem 3.1 (i)]{LeonoriPetitta2007} it is possible to  show that,  for each $p\in (1, \bar p]$, there exists a solution $u_p$ to \eqref{PpF}, that is, $u_p$ satisfies 
\begin{equation}
\label{ws}
\int_{\Omega}|\nabla u_p|^{p-2}\nabla u_p\cdot \varphi\,dx=\int_{\Omega}|\nabla u_p|^{p-2}\nabla u_p\cdot {\bf F}\varphi\,dx+\int_{\Omega}f\varphi\,dx\quad\hbox{ for all }\;\varphi\in W^{1,p}_0(\Omega).
\end{equation}
For the sake of completeness we sketch the proof of this result in  Appendix \ref{App}.

Moreover, the family $(u_p)_{p>1}$ is uniformly bounded in $p$. Indeed, use $u_p$ as test function in \eqref{ws} to obtain 
\begin{equation}
\label{upF}
\int_{\Omega}|\nabla u_p|^p\,dx=\int_{\Omega}|\nabla u_p|^{p-2}\nabla u_p\cdot {\bf F}u_p\,dx+\int_{\Omega}fu_p\,dx.
\end{equation}
From Proposition \ref{Inq}, H\"older inequality and Proposition \ref{SL}, it follows that
\begin{align*}
\label{Est0}
&\int_{\Omega}|\nabla u_p|^{p-2}\nabla u_p\cdot {\bf F}u_p \,dx\leq \int_{0}^{\infty}(|\nabla u_p|^{p-1})^*(t)|{\bf F}|^*(t)|u_p|^*(t)\,dt\nonumber\\
&\leq \||{\bf F}|\|_{L^{N,\infty}(\Omega)}\left(\int_{0}^{\infty}[(|\nabla u_p|^{p-1})^*(t)]^{p/(p-1)}\,dt\right)^{(p-1)/p}\left(\int_{0}^{\infty}[t^{1/{p^*}}|u_p|^*(t)]^p\,\frac{dt}{t}\right)^{1/p}\\
\nonumber&\leq \gamma_p\||{\bf F}|\|_{L^{N,\infty}(\Omega)}\int_{\Omega}|\nabla u_p|^p\,dx,
\end{align*}
 where $\gamma_p=p/[(N-p)C_N^{1/N}]$.
Thus, using \eqref{HFf}, \eqref{upF} becomes
\begin{equation}
\label{Est1}
\int_{\Omega}|\nabla u_p|^p\,dx\leq \left(\frac{S_N\|f\|_{L^{N}(\Omega)}}{1-\gamma_p\||{\bf F}|\|_{L^{N,\infty}(\Omega)}}\right)^{p/(p-1)}|\Omega|\leq \left(\frac{1-\gamma_1\||{\bf F}\|_{L^{N,\infty}(\Omega)}}{1-\gamma_p\||{\bf F}|\|_{L^{N, \infty}(\Omega)}}\right)^{p/(p-1)}|\Omega|.
\end{equation}
Since
$$
\lim_{p\to 1^+}\left(\frac{1-\gamma_1\||{\bf F}\|_{L^{N,\infty}(\Omega)}}{1-\gamma_p\||{\bf F}|\|_{L^{N, \infty}(\Omega)}}\right)^{p/(p-1)}=e^{(N\gamma_1\||{\bf F}|\|_{L^{N,\infty}(\Omega)})/[(N-1)(1-\gamma_1\||{\bf F}|\|_{L^{N,\infty}(\Omega)})]},
$$
there exists $\bar p>1$ such that
\begin{equation*}
\label{CF}
\int_{\Omega}|\nabla u_p|^p\,dx\leq C\quad\hbox{ for all }\; p\in (1,\bar p],\quad\hbox{ for some }\; C>0.
\end{equation*}
But this, in turn, up  to subsequences, implies that there exists $u\in BV(\Omega)$ such that 
\begin{eqnarray}
\label{CvU}
\nonumber\nabla u_p\rightharpoonup Du&*-\hbox{ weakly in }\; \mathcal{M}(\Omega,\re^N)\\
u_p\to u & \hbox{ in }\; L^q(\Omega),\; \forall \, q\in [1, 1^*),\quad\hbox{ a.e. in }\;\Omega\\
\nonumber u_p\rightharpoonup u& \hbox{ in }\; L^{1^*}(\Omega),
\end{eqnarray}
as $p\to 1^+$.

\subsection{Proof of the existence result for a generic drift}

Here is our notion of solution of problem \eqref{P1F} that naturally extends Definition \ref{DefSol}. 
\begin{definition}
\label{DSF}
We say that $u\in BV(\Omega)$ is a solution to \eqref{P1F} if there exist ${\bf z}\in X_{\mathcal{M}}(\Omega)$ with $\|{\bf z}\|_{L^\infty(\Omega,\re^N)}\leq 1$   such that
\begin{itemize}
\item [$(i)$]  $ -\hbox{div}\,{\bf z}={\bf z}\cdot{\bf F} +f\quad\hbox{ in }\; \mathcal{D}'(\Omega)$;
\item [$(ii)$] $ ({\bf z}, DT_k(u))= |DT_k(u)|\quad \hbox{ as Radon measures on }\Omega\quad\hbox{ for all }\;k>0$;
\item [$(iii)$] $\left[{\bf z}, \nu\right]\in \hbox{Sgn}\,(-u)\quad \mathcal{H}^{N-1}-\hbox{ a.e. in }\;\partial\Omega$.
\end{itemize}
\end{definition}

Our main existence result  in presence of a generic drift term ${\bf F}$ is  the following.  
\begin{theorem}
\label{TF}
Let $f\neq 0$ in  $L^N(\Omega)$ and let ${\bf F}\colon \Omega\to\mathbb{R}^N$ be a measurable vectorial function such that $|{\bf F}|\in L^{N, \infty}(\Omega)$. If $f$ and ${\bf F}$ satisfy \eqref{HFf} then there exists a solution to \eqref{P1F} in the sense of Definition \ref{DSF}. 
\end{theorem}

In oder to prove the above theorem we need the following.

\begin{lemma}
\label{Z}
For each $p>1$, let $u_p$ be a solution of \eqref{PpF}. Then there exists $ {\bf z}\in L^\infty (\Omega,\re^N)$ with $ \|{\bf z}\|_{L^\infty(\Omega, \mathbb{R}^N)}\leq 1$ such that, up  to subsequences, for any $1<q<\infty$
$$
|\nabla u_p|^{p-2}\nabla u_p\rightharpoonup {\bf z}\quad\hbox{ in }\; L^q(\Omega)\quad\hbox{ as }\; p\to 1^+.
$$
\end{lemma}
\dem
Let $q>1$. Using H\"older inequality and \eqref{Est1}, we obtain that
\begin{equation}
\label{Est2}
\left(\int_{\Omega}||\nabla u_p|^{p-2}\nabla u|^q\,dx\right)^{1/q}\leq \left(\frac{1-\gamma_1\||{\bf F}|\|_{L^{N, \infty}(\Omega)}}{1-\gamma_p\||{\bf F}|\|_{L^{N, \infty}(\Omega)}}\right)|\Omega|^{1/q}.
\end{equation}
But this, in turn, implies that there exists ${\bf z}_q\in L^q(\Omega, \mathbb{R}^N)$ such that, up  to subsequences,
$$
|\nabla u_p|^{p-2}\nabla u_p\rightharpoonup {\bf z}_q\quad\hbox{ in }\; L^q(\Omega, \mathbb{R}^N)\quad\hbox{ as }\; p\to 1^+
$$
 and, by a diagonal argument, there exists ${\bf z}$ independent of $q$ such that
$$
|\nabla u_p|^{p-2}\nabla u_p\rightharpoonup {\bf z}\quad\hbox{ in }\; L^q(\Omega, \mathbb{R}^N)\quad\hbox{ as }\; p\to 1^+.
$$

Passing to limit in \eqref{Est2}, as $p\to 1^+$, we obtain that
$$
\|{\bf z}\|_{L^q(\Omega,\re^N)}\leq |\Omega|^{1/q}\quad\hbox{ for all }\; q>1.
$$
So, letting $q\to \infty$,
$$
\|{\bf z}\|_{ L^\infty(\Omega,\re^N)}\leq 1.
$$
\fim

\noindent {\it Proof of Theorm \ref{TF}} In order to show this theorem, we prove that $u$ and  ${\bf z}$ given found, respectively, in \eqref{CvU} and in  Lemma \ref{Z}  satisfy $(i)-(iii)$ of Definition \ref{DSF}.

\noindent{\it Proof of  $(i)$.} It result by taking the  limit in \eqref{ws} as $p\to 1^+$, and using  Lemma \ref{Z}. In fact, as $|{\bf F}|\in L^{N,\infty}(\Omega)\subset L^{q'}(\Omega)$ and $|\nabla u_p|^{p-2}\nabla u_p\rightharpoonup {\bf z}$ in $ L^q(\Omega, \mathbb{R}^N)$,  one has that
$$
\lim_{p\to1^+}\int_{\Omega}|\nabla u_p|^{p-2}\nabla u_p\cdot {\bf F}\varphi\,dx=\int_{\Omega}{\bf z}\cdot {\bf F}\varphi\,dx,
$$
for all $\varphi \in C^\infty_c(\Omega)$. {Consequently, letting $p\to 1^+$ in \eqref{ws}, 
$$
\int_{\Omega}{\bf z}\cdot \nabla \varphi\,dx=\lambda\int_{\Omega}{\bf z}\cdot {\bf F}\varphi\,dx+\int_{\Omega}f\varphi\,dx\quad\hbox{ for all }\;\varphi\in C^\infty_c(\Omega),
$$
or equivalently
$$
-\hbox{div}\,{\bf z}=\lambda{\bf z}\cdot{\bf F}+f\quad\hbox{ in }\; \mathcal{D}'(\Omega).
$$
Note that, $-\hbox{div}\,{\bf z}\in L^{N,\infty}(\Omega)\subset \mathcal{M}(\Omega)$ and so, Anzellotti's theory is available.
}

\noindent{\it Proof of  $(ii)$.}  {Use \eqref{Dist} to obtain}
$$
({\bf z}, DT_k(u))\leq |DT_k(u)|\quad\hbox{ as Radon measures on }\;\Omega\quad\hbox{ for all }\; k>0.
$$
Let $0\leq \varphi\in C^\infty_c(\Omega)$. In order to proof the reverse inequality, we  use $\varphi T_k(u_p)$ as test function in \eqref{ws} to obtain
\begin{align*}
&\int_{\Omega}\varphi|\nabla T_k(u_p)|^p\,dx +\int_{\Omega}T_k(u_p)|\nabla u_p|^{p-2}\nabla u_p\cdot\nabla \varphi\,dx \\ &=\int_{\Omega}|\nabla u_p|^{p-2}\nabla u_p\cdot {\bf F}\varphi T_k(u_p)\,dx+\int_{\Omega}f \varphi T_k(u_p)\,dx.
\end{align*}
Now we let  $p\to 1^+$; using the weak lower semi-continuity, Young's inequality, and  Lemma \ref{Z},   one gets 
\begin{align*}
\int_{\Omega}\varphi|DT_k(u)|&\leq \liminf_{p\to 1^+}\int_{\Omega}\varphi|\nabla T_k(u_p)|\,dx\\
&\leq \liminf_{p\to 1^+}\left(\frac{1}{p}\int_{\Omega}\varphi|\nabla T_k(u_p)|^p\,dx+\frac{p-1}{p}\int_{\Omega}\varphi\,dx\right)\\
&= \lim_{p\to 1^+}\left(\frac{1}{p}\int_{\Omega}|\nabla u_p|^{p-2}\nabla u_p\cdot {\bf F}\varphi T_k(u_p)\,dx+\frac{1}{p}\int_{\Omega}f\varphi T_k(u_p)\,dx\right.\\
&\quad\left.-\frac{1}{p}\int_{\Omega}T_k(u_p)|\nabla u_p|^{p-2}\nabla u_p\cdot \nabla\varphi\,dx+\frac{p-1}{p}\int_{\Omega}\varphi\,dx\right)\\
&\leq \int_{\Omega}{\bf z}\cdot{\bf F}\varphi T_k(u)\,dx+\int_{\Omega}f\varphi T_k(u)\,dx-\int_{\Omega} T_k(u){\bf z}\cdot \nabla \varphi\,dx\\
&=-\int_{\Omega}\varphi T_k(u)\hbox{div}\,{\bf z}\,dx-\int_{\Omega}T_k(u){\bf z}\cdot\nabla \varphi\,dx\\
&=\int_\Omega ({\bf z}, DT_k(u))\varphi.
\end{align*}
Hence the contrary inequality holds.

\noindent{\it Proof of  $(iii)$.} Observe that, by Proposition \ref{Prop2}, we have that
\begin{equation}\label{reverse}
\int_{\partial\Omega}|T_k(u)|\,d\mathcal{H}^{N-1}+\int_{\partial\Omega}[{\bf z}, \nu] T_k(u)\,d\mathcal{H}^{N-1}\geq 0. 
\end{equation}
In order to show the reverse inequality, use $T_k(u_p)$ as test function in \eqref{ws} to obtain 
$$
\int_{\Omega}|\nabla T_k(u_p)|^p\,dx=\int_{\Omega}|\nabla u_p|^{p-2}\nabla u_p\cdot {\bf F} T_k(u_p)\,dx+\int_{\Omega}fT_k(u_p)\,dx.
$$
By Young's inequality and the fact that $u_p=0$ on $\partial\Omega$, it follows that
\begin{align*}
&\int_{\Omega}|\nabla T_k(u_p)|\,dx +\int_{\partial\Omega}|T_k(u_p)|\,d\mathcal{H}^{N-1} \\&\leq \frac{1}{p}\int_{\Omega}|\nabla u_p|^{p-2}\nabla u_p \cdot {\bf F} T_k(u_p)\,dx+\frac{1}{p}\int_{\Omega}f T_k(u_p)\,dx+\frac{p-1}{p}|\Omega|.
\end{align*}
Letting $p\to 1^+$, by lower semi-continuity, we obtain that
\begin{align*}
\int_{\Omega}|DT_k(u)|+\int_{\partial\Omega}|T_k(u)|\,d\mathcal{H}^{N-1}&\leq \int_{\Omega} {\bf z}\cdot{\bf F} T_k(u)\,dx +\int_{\Omega}f T_k(u)\,dx\\
&=-\int_{\Omega} T_k(u)\hbox{div}\,{\bf z}\,dx\\
&=\int_{\Omega}({\bf z}, DT_k(u))-\int_{\Omega}[{\bf z}, \nu] T_k(u)\,d\mathcal{H}^{N-1}.
\end{align*}
But this, by $(ii)$, implies that the reverse  to \eqref{reverse} holds,  and so
\begin{equation}
\label{Est4}
\int_{\partial\Omega}|T_k(u)|\,d\mathcal{H}^{N-1}+\int_{\partial\Omega}[{\bf z}, \nu] T_k(u)\,d\mathcal{H}^{N-1}=0.
\end{equation}
As $T_k(u)\to u$ in sense of strict convergence, by continuity of the trace operator in $BV(\Omega)$, we have that  
$$
T_k(u)\to u \quad\hbox{ in }\; L^1(\partial\Omega)\quad\hbox{ as }\; k\to \infty.
$$
Hence, letting $k\to \infty$ in \eqref{Est4},
$$
\int_{\partial\Omega}|u|\,d\mathcal{H}^{N-1}+\int_{\partial\Omega}[{\bf z}, \nu] u\,d\mathcal{H}^{N-1}=0.
$$
and so, $(iii)$ holds.

Therefore $(i)-(iii)$ of Definition \ref{DSF} hold and so, $u$ is a solution solution to \eqref{P1F}.
\fim

\begin{remark}\label{remextend}
 If one assume that $f\in L^q(\Omega)$ with $q>N$, then  by Stampacchia's argument, as in the proof of Theorem \ref{bo}, one can show that $(u_p)_{p>1}$ is uniformly bounded and that $u\in L^\infty(\Omega)$.
 Also observe that, if  $f\in L^{N, \infty}(\Omega)$ is such that
\begin{equation*}
(\||{\bf F}|\|_{L^{N,\infty}(\Omega)}+\|f\|_{L^{N,\infty}(\Omega)})\leq \frac{1}{\gamma_1}.
\end{equation*}
By combining the very same arguments of the proofs of Theorem \ref{exiinf} and  of Theorem \ref{TF}, one can show that there exists ${\bf z}\in X_{\mathcal{M}}(\Omega)$ with $\|{\bf z}\|_{L^\infty(\Omega,\mathbb{R}^N)}\leq 1$  such that
\begin{itemize}
\item [$(i)$]  $-\hbox{div}\,{\bf z}= {\bf z}\cdot{\bf F}+f \quad\hbox{ in }\; \Omega$;
\item [$(ii)$]  $({\bf z}, DT_k(u))=|DT_k(u)|\quad\hbox{ as Radon measures on }\; \Omega, \;\hbox{ for all }\;k>0$;
\item [$(iii)$] $[{\bf z}, \nu]\in \hbox{Sgn}(-u)\quad\mathcal{H}^{N-1}-\hbox{a.e. on }\; \partial\Omega$,
\end{itemize}
and that the analogous of Theorem \ref{AsympB} also holds. 
\end{remark}

\section{Non-existence and optimality. Some explicit examples}\label{sec7} 
In this section, we give some explicit examples provided  assumption \eqref{HL} is either satisfied or not. We start with the following example, in which  the equality sign  in \eqref{HL} holds; in this case one can construct an example in which the limit $u$ of $(u_p)_{p>1}$ found in Theorem \ref{AsympB} is a  nontrivial solution of \eqref{P}.

\begin{example} Let $R>0$ and $\lambda=-(N-2)$. Consider the problem 
\begin{equation}
\label{TP}
\left\{
\begin{array}{rclr}
-\Delta_p u&=&\lambda |\nabla u|^{p-2}\nabla u\cdot\frac{ x}{|x|^2}+\frac{1}{|x|}&\hbox{ in }\; B_R(0)\\
u&=&0&\hbox{ on }\;\partial B_R(0).
\end{array}
\right.
\end{equation}
One can easily show, that the function
$$
u_p(x)=\frac{R^{N-1}}{N-1}\left[1-\left(\frac{|x|}{R}\right)^{N-1}\right]
$$
is a solution of \eqref{TP} and
$$
\frac{|\lambda|}{N-1}+\gamma_1\left\|\frac{1}{|x|}\right\|_{L^{N, \infty}(\Omega)}=1,
$$
where $\gamma_1$ is as in \eqref{Tg}. 

Clearly,
$$
\lim_{p\to 1^+}u_p(x)=\frac{R^{N-1}}{N-1}\left[1-\left(\frac{|x|}{R}\right)^{N-1}\right]=u(x) \quad \hbox{ for all } \; x\in B_R(0), $$ and  $$ {\bf z}=-\frac{x}{|x|}\quad\hbox{ for all }\;x\in B_R(0),
$$
satisfy the following conditions 
\begin{eqnarray*}
-\hbox{div}\,{\bf z}&=&\lambda{\bf z}\cdot\frac{x}{|x|^2}+\frac{1}{|x|}\qquad\hbox{ in }\; B_R(0)\\
({\bf z}, Du)&=&|Du|\qquad\hbox{ as  Radon measures in }\; B_R(0)\\
\left[{\bf z}, \frac{x}{|x|}\right]&\in& \hbox{Sgn}(-u)\qquad\mathcal{H}^{N-1}-\hbox{a.e.}\;x\in \partial B_R(0),
\end{eqnarray*}
that is, $u$ is a solution to problem
\begin{equation*}
\left\{
\begin{array}{rclc}
-\Delta_1u&=&\lambda \frac{Du}{|Du|}\cdot \frac{x}{|x|^2}+\frac{1}{|x|}&\hbox{ in }\; B_R(0)\\
              u&=&0&\hbox{ on }\; \partial B_R(0).
\end{array}
\right.
\end{equation*}
\end{example}

On the other hand, the existence of a non-trivial solution of problem \eqref{P} when the equality  in \eqref{HL} holds   is not always true, as shown in the next example.

\begin{example}
For $N>3$, let $0<\beta < e^{(\left(N/(N-1)\right)^2-1)C_N}$. Let $f\colon B_R(0)\to \mathbb{R}$ be a function defined by
\begin{equation*}
f(x)=\left\{
\begin{array}{lcr}
\frac{1}{R/\beta}&\hbox{if}& |x|\leq \frac{R}{\beta}\\
\frac{1}{|x|}&\hbox{if}& \frac{R}{\beta}\leq |x|\leq R.
\end{array}
\right.
\end{equation*}
Note that, for all $s>0$, the distribution function of $f$ can be write as
\begin{equation*}
\alpha_f(s)=\left\{
\begin{array}{lcr}
C_N R^N&\hbox{ if } & 0<s<\frac{1}{R}\\
\frac{C_N}{s^N}&\hbox{ if } & \frac{1}{R}\leq s< \frac{\beta}{R}\\
C_N\left[\frac{1}{s^N}-\left(\frac{R}{\beta}\right)^N\right]& \hbox{if }& s\geq \frac{\beta}{R}.
\end{array}
\right.
\end{equation*}
The non-increasing rearrangement of $f$ is given by
$$
f^*(t)=\sup\{ s>0\colon \alpha_f(s)>t\}=\frac{C^{1/N}_N}{t^{1/N}} \quad\hbox{ for all }\;t\in (0, C_NR^N),
$$
and so, 
$$
 \|f\|_{L^{N, \infty}(\Omega)}=C_N^{1/N}.
$$
Thus, for $\lambda =-(N-2)$, 
$$
\frac{|\lambda|}{N-1}+\gamma_1\|f\|_{L^{N, \infty}(\Omega)}=1,
$$
where $\gamma_1=[(N-1)C^{1/N}_N]^{-1}$, that is,  \eqref{HL} is satisfied with the equality sign. But, as $\beta< e^{((N/(N-1))^N-1)C_N}$,
$$
\frac{|\lambda|}{N-1}+S_N\|f\|_{L^N(B_R(0))}\leq \frac{N-2}{N-1}+\frac{1}{NC_N^{1/N}}(C_N+\ln(\beta))^{1/N}< 1
$$
 and so, by Theorem \ref{AsympB}
$$
\lim_{p\to 1^+}u_p(x)=0\quad\hbox{ a.e. }\; x\in B_R(0).
$$
\end{example}

The following example shows the optimality of the condition \eqref{HL}.

\begin{example}
Let $\lambda=-(N-2)$. Consider the datum $f(x)=\alpha/|x|$, with $\alpha>0$. Note that
$$
u_p(x)=\frac{\alpha^{1/(p-1)}R^{N-1}}{N-1}\left[1-\left(\frac{|x|}{R}\right)^{N-1}\right]\quad\hbox{ for all }\; p>1.
$$ 
is a solution to
\begin{equation*}
\left\{
\begin{array}{rclr}
-\Delta_p u&=&\lambda |\nabla u|^{p-2}\nabla u\cdot\frac{x}{|x|^2}+\frac{\alpha}{|x|}&\hbox{ in }\; B_R(0)\\
u&=&0&\hbox{ on }\;\partial B_R(0).
\end{array}
\right.
\end{equation*}
However
\begin{equation*}\scriptsize{
\lim_{p\to 1^+}u_p(x)=\left\{
\begin{array}{lccr}
0&\hbox{ if }& \alpha <1&\iff |\lambda|/(N-1)+\gamma_1\|\alpha/|x|\|_{L^{N,\infty}(\Omega)}<1\\
R^{N-1}/[(N-1)(1-(|x|/R)^{N-1})]& \hbox{ if }& \alpha =1& \iff |\lambda|/(N-1)+\gamma_1\|\alpha/|x|\|_{L^{N,\infty}(\Omega)}=1\\
+\infty& \hbox{ if }& \alpha >1&\iff |\lambda|/(N-1)+\gamma_1\|\alpha/|x|\|_{L^{N,\infty}(\Omega)}>1.
\end{array}
\right.}
\end{equation*}
\end{example}

\subsection{The case $1-N<\lambda<0$ and an associated singular weighted problem}\label{sbs6}

Let us conclude by  focus our attention to   the particular  case of a negative parameter $\lambda>1-N$. In this case in fact, problem \eqref{P}  exhibits  an intriguing connection with some singular elliptic problems that are naturally studied in weighted Sobolev space. More precisely, if we multiply by $|x|^{\lambda}$ in \eqref{P}, we obtain that
$$
-\hbox{div}\,\left(|x|^\lambda\frac{Du}{|Du|}\right)=|x|^\lambda f\quad\hbox{ in }\;\Omega.
$$
Denote $a=-\lambda$, so $ 0<a<N-1$ and our   problem becomes
\begin{equation}
\label{Pa}
\left\{
\begin{array}{rclc}
-\hbox{div}\,\left(\frac{1}{|x|^a}\frac{Du}{|Du|}\right)&=&\frac{1}{|x|^a}f&\hbox{ in }\; \Omega\\
              u&=&0&\hbox{ on  }\;\partial\Omega.
\end{array}
\right.
\end{equation}
 
 Problems as is \eqref{Pa} are well-known in the literature (see for instance  \cite{OrtizPimentaSegura2021} and references therein) and they are well-posed  in  weighted Sobolev spaces. With the help of the extension of the the Anzellotti theory  developed in \cite{OrtizPimentaSegura2021}, in fact, one can   show the  existence of a solution to \eqref{Pa}. More specifically, it can be seen that there exists $u\in BV(\Omega)$ and ${\bf z}\in L^\infty(\Omega, \mathbb{R}^N)$ with $\|{\bf z}\|_{L^\infty(\Omega, \mathbb{R}^N)}\leq 1$ such that
\begin{equation*}
\left\{
\begin{array}{lr}
-\hbox{div}\,\left(\frac{1}{|x|^a}{\bf z}\right)=\frac{1}{|x|^a}f& \hbox{ in }\; \mathcal{D}'(\Omega)\\
\frac{1}{|x|^a}({\bf z}, Du)=\frac{1}{|x|^a}|Du|&\hbox{ as Radon measures on  }\; \Omega\\
\left[{\bf z}, \nu\right]\in \hbox{Sgn}(-u)&\mathcal{H}^{N-1}-\hbox{ a.e. in }\; \partial\Omega.
\end{array}
\right.
\end{equation*}

\appendix\section{Existence for the approximating problem }\label{App}

In  \cite[Theorem 3.1 (i)]{LeonoriPetitta2007}, the existence of a weak solution for problem \eqref{PpF} was proven under the slightly stronger assumption that 
$$
|{\bf F}(x)|\leq \frac{B}{|x|}\, \ \forall\ x\in\Omega\,.
$$
If one only assume $|{\bf F}|\in L^{N, \infty}(\Omega)$ {such that
\begin{equation}
\label{FGp} 
\||{\bf F}|\|_{L^{N,\infty}(\Omega)}< \frac{1}{\gamma_p} \quad\hbox{ with }\; \gamma_p=\frac{p}{(N-p)C^{1/N}_N},
\end{equation}}
the proof is very similar; although, for the sake of completeness here we sketch it by emphasizing the differences relying in passing to the limit in the lower order term.   

\begin{theorem}
\label{TFp}
{Let $1<p<N$}. Let {$f\in L^m(\Omega), \; m>N/p$},  and let ${\bf F}\colon \Omega\to\mathbb{R}^N$ be a measurable vectorial function such that $|{\bf F}|\in L^{N, \infty}(\Omega)$ {satisfying  \eqref{FGp}}. Then there exists a weak  solution ${u}\in W^{1,p}(\Omega)$ of \eqref{PpF} in the sense of \eqref{ws}. 
\end{theorem}

\dem
For $f_n=T_n(f)$,  consider the problem 
{
\begin{equation*}
\left\{
\begin{array}{rclc}
-\Delta_p u_n &=& \frac{|\nabla u_n |^{p-2}\nabla u_n \cdot {\bf F}}{1+1/n|\nabla u_n|^{p-1}|{\bf F}|}+ f_n& \hbox{ in }\; \Omega\\
u_n&=&0& \hbox{ on }\;\partial\Omega.
\end{array}
\right.
\end{equation*}
}
By classical Leray-Lions theory, it well known that the above problem has solution  $u_n\in {W^{1,p}_0}(\Omega)$ satisfying  
\begin{equation}
\label{wsn}
\int_{\Omega}|\nabla u_n|^{p-2}\nabla u_n\cdot \nabla \varphi\,dx=\int_{\Omega}\frac{|\nabla u_n|^{p-2}\nabla u_n\cdot {\bf F}}{1+1/n|\nabla u_n|^{p-1}|{\bf F}|}\varphi\,dx+\int_{\Omega}f_n\varphi\,dx{\quad\hbox{ for all }\; \varphi\in W^{1,p}_0(\Omega)}.
\end{equation}

 Use $u_n$ as test function above to obtain, {by H\"older and Sobolev inequalities,  and \eqref{FGp},}
 
\begin{equation*}
\|u_n\|_{W^{1,p}_0(\Omega)}\leq{\left(\frac{S_{N,p}\|f\|_{L^m(\Omega)}}{1-\gamma_p\||{\bf F}\|_{{L^{N, \infty}(\Omega)}}}\right)^{{1/(p-1)}}|\Omega|^{(1-1/m-1/{p^*})/(p-1)}}\quad\hbox{ for all }\;n\in \mathbb{N}.
\end{equation*}

Consequently, there exist $u\in W^{1,p}_0(\Omega)$ such that, up to subsequences,
\begin{align*}
u_n\rightharpoonup u\quad\hbox{in }\; W^{1,p}_0(\Omega)\\
u_n\to u\quad\hbox{ in }\; L^q(\Omega)\quad\hbox{ for all }\; q\in [p,p^*),
\end{align*}
{as $n\to \infty$}. But this, in turn, allows to apply classical a.e. convergence of the gradients results (see \cite{BoccardoMurat1992}), i.e. 
$$
\nabla u_n\to \nabla u\quad\hbox{a.e. in }\; \Omega.
$$
Thus, by Vitali theorem, 
\begin{equation}
\label{Cq}
\nabla u_n\to \nabla u\quad\hbox{ in }\; (L^q(\Omega))^N\quad\hbox{ for all }\; q<p.
\end{equation}
which, in turn, implies that there exists $g\in L^q(\Omega)$ such that
$$
|\nabla u_n|\leq g\quad\hbox{a.e. in }\; \Omega.
$$

In order to pass to the limit in \eqref{wsn} as $n\to \infty$, we consider $N/(N-1)<q< p/(p-1)$ and so, $p<q'<N$. Note that $|{\bf F}|\in L^{N,\infty}(\Omega)\subset L^{q'}(\Omega)$. Thus
$$
|\nabla u_n|^{p-1}|{\bf F}|\leq g(x)^{p-1}|{\bf F}|\leq \frac{1}{q}g^{q(p-1)}+\frac{1}{q'}|{\bf F}|^{q'}\in L^1(\Omega).
$$
By Lebesgue dominated convergence theorem, it follows that
\begin{equation}
\label{CF1}
\lim_{n\to \infty}\int_{\Omega}\frac{|\nabla u_n|^{p-2}\nabla u_n\cdot {\bf F}}{1+1/n|\nabla u_n|^{p-1}|{\bf F}|}\varphi\,dx=\int_{\Omega}|\nabla u|^{p-2}\nabla u\cdot {\bf F}\varphi\,dx.
\end{equation}

Letting $n\to \infty$ in \eqref{wsn}, from \eqref{Cq} and \eqref{CF1}, we obtain
$$
\int_{\Omega}|\nabla u|^{p-2}\nabla u\cdot \nabla \varphi\,dx=\int_{\Omega}|\nabla u|^{p-2}\nabla u\cdot {\bf F}\varphi\,dx+\int_{\Omega}f\varphi\,dx\quad\hbox{ for all }\; \varphi\in W^{1,p}_0(\Omega),
$$
and this concludes the proof. 
\fim

\bk 
{\bf Acknowledgment.} Francesco Petitta is partially supported by  the Gruppo Nazionale per l’Analisi Matematica, la Probabilità e le loro Applicazioni (GNAMPA) of the Istituto Nazionale di Alta Matematica (INdAM). 
Juan C. Ortiz Chata is supported by FAPESP 2021/08272-6 and 2022/06050-9, Brazil.


\begin{thebibliography}{99}

\bibitem{Alvino1977}{\sc Alvino, A.:} {\em Sulla diseguaglianza di {S}obolev in spazi di {L}orentz},  Bolletino dell'Unione Matematica Italiana, 14-A(5) (1977) no. 1, 148--156.

\bibitem{AmbrosioFuscoPallara}{\sc Ambrosio, L., Fusco, N., Pallara, D.:} {Functions of bounded variation and free discontinuity problems}, Oxford University Press, Oxford, 2000.

\bibitem{acm2001}{\sc Andreu, F., Ballester, C., Caselles, V., Maz\'on, J.M.:} {\em The {D}irechlet problem for the total variation flow}, J. Funct. Anal. 180 (2001), no. 2, 347--403.


\bibitem{AndreuCasellesMazon2004}{\sc Andreu, F., Caselles, V. and Maz\'on, J.M.:} {Parabolic quasilinear equations minimizing linear growth functionals}, {P}rogress in mathematics 223, {B}irkh\"auser  {V}elarg, {B}asel, 2004.


\bibitem{4} {\sc Andreu F.,  Dall'Aglio A.,  and  Segura de Le\'on S.}: {\em Bounded solutions to the $1$-Laplacian equation with a critical gradient term}, Asymptotic Anal., 80 (2012), no. 1-2, 21--43.

\bibitem{AndreuMazonMollCaselles2004}{\sc Andreu, F.; Maz\'on, J.M.; Moll, S.; Caselles, V.:} {\em The minimizing total variation flow with measure initial conditions}, Comm. Contemp. Math. 6 (2004), no. 3, 431--494.

\bibitem{Anzellotti1983} {\sc Anzellotti, G.:} {\em Pairings between measures and bounded functions and compensated compactness}, Ann. Mat. Pura Appl., 135 (1983),  no. 1, 293--318.

\bibitem{AttouchButtazzoMichaille2006} {\sc Attouch, H.; Buttazzo, G.; Michaille, G.:} {Variational analysis in Sobolev spaces: applications to PDEs and optimization.} MPS-SIAM, Philadelphia, 2006. 

\bibitem{bop} {\sc  Balducci F., Oliva F.,  Petitta F.:} {\em  Finite energy solutions for nonlinear elliptic equations with competing gradient, singular and $L^1$ terms}, J. Differential Equations,  391 (2024), 334--369.


\bibitem{BCRS} {\sc  Bertalmio M.; Caselles V.; Roug\'e B.;   Sol\'e A. }: 
{\em TV based image restoration with local constraints.}
Special issue in honor of the sixtieth birthday of Stanley Osher, 
J. Sci. Comput. 19 (2003), No. 1-3, 95--122.

\bibitem{BMMP2003}{\sc Betta, M.; Mercaldo, A.; Murat, F.; Porzio, M.M.:} {\em Existence of renormalized solutions to non equations with a lower-order term and righ measure}, J. Math. Pures Appl. (9), 82 (2003), no. 1, 90--124.

\bibitem{BG1989}{\sc Boccardo, L.; Giachetti, D.:} {\em Existence results via regularity for some nonlinear elliptic problems}, Comm. Partial Differential Equations, 14 (1989), no. 5, 663--680.

\bibitem{BoccardoMurat1992}{\sc Boccardo, L.; Murat, F.:} {\em Almost everywhere convergence of the gradients of solutions to elliptic and parabolic equations}, Nonlinear Anal. 19 (1992), no. 19, 581--597.

\bibitem{BoccardoOrsinaPeral2006}{\sc Boccardo, L; Orsina, L; Peral, I.:} { \em A remark on existence and optimal summability of solutions of elliptic problems involving Hardy potential}, Discrete Contin. Dyn Syst. 16 (2006), 513--523.

\bibitem{CKN1984}{\sc Caffarelli, L.; Kohn, R.; Nirenberg, L.:} {\em First order interpolation inequalities with weights}, Compositio Math. 53 (1984), no. 3, 259--275.

\bibitem{Caselles2011}{\sc Caselles, V.:} {\em On the entropy conditions for some flux limited equations}, J. Differential Equations, 250 (2011), no. 8, 3311--3348.

\bibitem{crt}  {\sc Cassani D.,  Ruf B., and Tarsi C.:} {\em Equivalent and attained version of Hardy’s
inequality in $\mathbb{R}^N$,} J. Funct. Anal., 275 (2018), no. 12, 3303--3324. 

\bibitem{CF} {\sc Chen, G.-Q; Frid, H.:} {\em Divergence-measure fields and hyperbolic conservation laws}, Arch. Ration. Mech. Anal. 147 (1999), no. 2, 89--118.

\bibitem{ct}{\sc Cicalese, M.; Trombetti, C.:} {\em Asymptotic behaviour of solutions to {$p$}-{L}aplacian equation}, Asympt. Anal. 35 (2003), no. 1, 27--40.

\bibitem{DAnconaFanelli2007}{\sc D'Ancona, P.; Fanelli, L.:} {\em Decay Estimates for the Wave and Dirac Equations with a Magnetic Potential}, Comm. Pure Appl. Math, 60 (2007), no. 3, 357--392.

\bibitem{dgop} {\sc De Cicco V., Giachetti, D., Oliva, F., Petitta, F.: }{\em The Dirichlet problem for singular elliptic
equations with general nonlinearities}, Calc. Var. Partial Differential Equations, 58 (2019), no. 4, 1--40.

\bibitem{DV1995}{\sc Del Vecchio, T.:} {\em Nonlinear elliptic equations with measure}, Potential Anal. 4 (1995), no. 2, 185--203.


\bibitem{Demengel1999} {\sc   Demengel, F.:} {\em On some nonlinear partial differential equations Involving the {{\(``1"\)}}-{Laplacian} and critical {Sobolev} exponent}, ESAIM, Control Optim. Calc. Var., 4  (1999), 667-686.
,
\bibitem{EvansGariepy2015}{\sc Evans, L.C., Gariepy, R.F.:} {Measure theory and fine properties of functions}, Textbook in Matematics,  CRC Press, Boca Raton, FL (2015), xiv+299.

\bibitem{AzoreroAlonso1998}{\sc Garc\'{\i}a Azorero, J. P. and Peral Alonso, I.:} {\em Hardy inequalities and some critical elliptic and parabolic problems}. J. Differential Equations, 144 (1998), No. 2, 441--476.

\bibitem{gop} {\sc Giachetti D., Oliva F., and Petitta  F.}: {\em 1-{Laplacian} type problems with strongly singular nonlinearities and gradient terms}, Commun. Contemp. Math. 24 (2022), no.10 (10), Paper No. 2150081, 40 pp.

{\bibitem{HLP1964}{\sc Hardy, G.H.; Littlewood; P\'olya, G.:}{ \em Inequalities}, Cambridge Univ. Press, 1964.}\bk

\bibitem{Kawohl1990}{\sc Kawohl, B.:} {\em On a family of torsional creep problems}, J. Reine Angew. Math. 410 (1990),  1--22.

\bibitem{Kawohl1991}{\sc Kawohl, B.:} {\em From {$p-${L}aplace to mean curvature operator and related questions}. Progress in partial differential equations: the {M}etz surveys}, 249 (1991), 40--56.

\bibitem{lops}{\sc Latorre, M., Oliva F., Petitta F., Segura de Le\'{o}n, S.:} {\em The Dirichlet problem for the $1$-Laplacian with a general singular term and  L1 -data}
Nonlinearity, 34 (2021), no. 3, 1791--1816.
\bibitem{LMPP2011}{\sc Leonori, T.; Mart\'{\i}nez-Aparicio, P.J.; Primo, A.:} {\em Nonlinear elliptic equations with {H}ardy potential and lower order term with natural growth}, Nonlinear Anal. 74 (2011), no. 11, 3556--3569.

\bibitem{LeonoriPetitta2007}{\sc Leonori, T., Petitta, F.:} {\em Existence and regularity results for some singular elliptic problems}, Adv. Nonlinear Stud. 7 (2007), no. 3, 329--344.

\bibitem{MS2013}{\sc Maz\'{o}n, J.M., Segura de Le\'{o}n, S.:} {\em The {D}irichlet problem for a singular elliptic equation arising in the level set formulation of the inverse mean curvature flow}, Adv. Calc. Var. 6 (2013), no. 2, 123--164.

{ \bibitem{mst}{\sc Mercaldo, A., Segura de Le\'on, S., Trombetti, C.:} {\em On the behaviour of the solutions to $p$-Laplacian equations as $p$ goes to 1}. Publ. Mat. 52 (2008), no. 2, 377--411.\bk }
\bibitem{MST2009}{\sc Mercaldo, A., Segura de Le\'{o}n, S., Trombetti, C.:}{ \em On the solutions to 1-{L}aplacian equation with {$L^1$} data}, J. Funct. Anal. 256 (2009), no. 8, 2387--2416.

\bibitem{M} {\sc  Meyer Y.:}  {Oscillating Patterns in Image Processing and Nonlinear Evolution Equations: The Fifteenth Dean Jacqueline B. Lewis Memorial Lectures}.   Providence, RI: American Mathematical Society, 2001.

\bibitem{Miranda1974/75}{\sc Miranda, M.:} {\em Direchlet problem with {$L\sp{1}$} data  for the non-homogeneous minimal surface equation}, Indiana Univ. Math. J., 24, (1974/75), 227--241.

\bibitem{ONeil1968}{\sc O'Neil, R.:} {\em Integral transforms and tensor products on {O}rlicz spaces and $L(p, q)$ spaces}, J. Anal. Math. 21 (1968), 1--276.

\bibitem{OrtizHiroshi2024}{\sc Ortiz Chata, J.C.; Hiroshi Miyagaki, O.:}{ \em Existence of solutions to a $1$-Laplacian problem in $\mathbb{R}^N$ with unbounded weight on the nonlinearity}, preprint.

\bibitem{OrtizPetitta2024}{\sc Ortiz Chata, J. C., Petitta, F.:} {\em  Existence, non-existence and degeneracy of limit solutions to Laplace problems involving Hardy potentials as $p\to1^+$,} SIAM J. Math. Anal., in press (2024).

\bibitem{OrtizPimentaSegura2021}{\sc Ortiz Chata, J. C.; Pimenta, M.T.O.; Segura de Le\'on, S.:} {\em Anisotropic 1-{L}aplacian problems with unbounded weights}, NoDEA Nonlinear Differential Equations Appl., 28 (2021), No. 57, 40.
\bibitem{OsSe}
{\sc Osher S.; Sethian J.:} {\em Fronts propagating with curvature-dependent speed: algorithms based on Hamilton-Jacobi formulations}, Journal of Computational Physics, 79 (1998), no.1, 12--49. 

\bibitem{pe}{\sc Peetre, J.:}  {\em Espaces d'interpolation et th\'eor\`eme de Soboleff}, Ann. Inst. Fourier (Grenoble) 16 (1966), 279--317.

 \bibitem{ROF}{\sc   Rudin, L.I. Osher, S.; Fatemi, E.}: {\it Nonlinear total variation based noise removal algorithms}, 
   Physica D. 60 (1992), 259--268.

    \bibitem{Sapiro}  {\sc  Sapiro G.:}  {Geometric partial differential equations and image analysis,} Cambridge University Press, 2001.

\bibitem{Stampacchia1965}{\sc Stampacchia, G.:} {\em Le probl\'eme de {D}irichlet pour les \'equations elliptiques du second ordre \`a coefficients discontinus}, Ann. Inst. Fourier(Grenoble), 15  (1965), 189--258.

\bibitem{Talenti}{\sc Talenti, G.:} {\em Best constant in {S}obolev inequality}, Ann. Mat. Pura Appl. 110 (1976), 353--372.





\end{thebibliography}
\end{document}